\documentclass[10pt,a4paper]{amsart}

\usepackage[T1]{fontenc}
\usepackage[utf8]{inputenc}
\usepackage[normalem]{ulem}

\usepackage{amsmath}
\usepackage{amsthm}
\usepackage{mathrsfs}

\usepackage[margin=2cm]{geometry}
\usepackage{xcolor}
\usepackage{float}
\usepackage{comment}
\usepackage{tikz}
\usepackage{tikz-cd}
\usetikzlibrary{arrows.meta}
\usetikzlibrary{positioning}

\usepackage[sort&compress,numbers,square]{natbib}
\usepackage{hyperref}
\hypersetup{
  unicode,
  colorlinks=true,
  linkcolor=violet,
  citecolor=green!50!black,
  urlcolor=violet,
  pdfauthor={Author One, Author Two, Author Three},
  pdftitle={Dimension filtration of Leavitt path K-algebras and intersection with Thompson's groups},
  pdfsubject={Filtered growth of Leavitt path algebras and Thompson-type groups},
  pdfkeywords={Leavitt path algebra, growth series, Thompson groups}
}
\def\dirloopedge#1(#2,#3)(#4,#5){
    \pgfpointxy{#2}{#3}
    \pgfgetlastxy{\ax}{\ay}
    \pgfmathsetmacro\outangle{atan2(#2,#3)}
    \pgfmathsetmacro\inangle{atan2(#4,#5)}

    \draw [-latex] (#1) -- +(#2,#3) to [out=\outangle,in=\inangle,looseness=5]  +(#4,#5) -- (Ro);
}

\theoremstyle{plain}
\newtheorem{theorem}{Theorem}[section]
\newtheorem{cor}[theorem]{Corollary}
\newtheorem{lemma}[theorem]{Lemma}

\newtheorem{prop}[theorem]{Proposition}

\theoremstyle{definition}
\newtheorem{definition}[theorem]{Definition}
\newtheorem{example}[theorem]{Example}
\newtheorem{remark}[theorem]{Remark}

\theoremstyle{plain}
\newtheorem{teoInt}{Theorem}

\renewcommand{\L}{\mathbb{L}}
\newcommand{\Z}{\mathbb{Z}}
\newcommand{\Q}{\mathbb{Q}}

\newcommand{\gr}{\mathrm{gr}}
\newcommand{\GK}{\mathrm{GK}}
\newcommand{\entr}{\mathrm{h}_\mathrm{alg}}
\newcommand{\TFG}{\mathbf{F}(\mathcal{G}_E)}

\title{The standard filtration of Leavitt path algebras and topological full groups.}

\author{Conchita Martínez-Pérez \and Sofía Sirón}


\begin{document}

\thanks{\noindent Both authors were partially supported by the Spanish Government PID2024-155800NB-C32 and by Gobierno de Arag{\'o}n and Feder Funds (grant code: E22-26R: ``{\'A}lgebra y Geometr{\'i}a'')}

\begin{abstract} We consider several problems related to the growth of Leavitt path algebras with respect to  the natural length filtration. We define the generating function of the series of dimensions of the successive quotients of this filtration and show that it is a rational function with poles which are reciprocal of the eigenvalues of the adjacency matrix, a fact which is related to recent results by Bock,  Gil Canto,  Mart\'in Barquero,  Mart\'in Gonz\'alez, Ruiz Campos and Sebandal on algebraic entropy of Leavitt path algebras. We also consider the filtration induced in the topological full group of the Renault-Deaconu groupoid of the defining graph and study the growth of the associated series. In the particular case of Thompson's group $V$, we find double exponential upper and lower bounds for this growth. We also prove related results for the group $F$ and for some submonoids. 
\end{abstract}


\maketitle

\section{Introduction}
Leavitt algebras  $\mathbb L_r$ of type $(1,r)$ were introduced by Leavitt in the 60's as canonical examples of rings without the invariant basis number property. Later, Ara, Moreno and Pardo, and independently Abrams and Aranda Pino  generalised the construction and defined the notion of Leavitt path algebra $L_K(E)$  associated to a graph $E$ (here, $K$ is a field). The defining graph $E$ determines the algebraic properties of the algebra $L_K(E)$ and in this sense, Leavitt path algebras provide a natural transfer between combinatorics and algebra. Moreover, they are part of a family of objects deeply related to each other; such as the graph $C^*$-algebra, the Renault-Deaconu groupoid $\mathcal{G}_E$ or the associated topological full group $\TFG$. Steinberg has defined an analogue of the group algebra for groupoids, the so called Steinberg algebra, which in the case of $\mathcal{G}_E$ is isomorphic to $L_K(E)$ \cite{OFST}.

Leavitt path algebras admit a natural grading over the integers where the generators associated to ordinary edges of $E$ are assigned degree 1 and generators associated to ghost edges are assigned degree -1. This grading has been extensively studied and many results on Leavitt path algebras have also a graded version, but has the disadvantage that can not be used to study  problems related to growth because the homogeneous components are typically of infinite dimension. Alternatively, it is possible to define a filtration $L_K(E)_n$, $n\geq 0$ of a Leavitt path algebra essentially by assigning degree 1 to generators associated to both ordinary and ghost edges, in this case the terms  $L_K(E)_n$ are finite-dimensional if the graph is finite and as a consequence, one can look at the growth of the series  $\dim_KL_K(E)_n$ of dimensions. In particular, a natural question  is whether the growth of this series can be described using combinatorial data of the graph. One interesting invariant is the Gelfand-Kirillov dimension, defined as 
$$\GK(L_K(E))=\limsup_{n\to\infty}\frac{\log(\dim_KL_K(E)_n)}{\log(n)}.$$
 In \cite{Zelmanov}, Alahmadi, Alsulami,  Jain and Zelmanov show that $\GK(L_K(E))$ is finite if and only if the graph $E$ has the property that cycles are disjoint from each other (see Theorem \ref{thm:Zelmanov}). This condition is rather restrictive and in \cite{entropy1}, \cite{entropy2},  Bock,  Gil Canto,  Mart\'in Barquero,  Mart\'in Gonz\'alez, Ruiz Campos and Sebandal,
consider a finer invariant, the algebraic entropy,  
$$\entr(L_K(E))=\limsup_{n\to\infty}\frac{\log(\dim_K(L_K(E)_n/L_K(E)_{n-1}))}{n}$$
and show $\entr(L_K(E))=\log(\rho(A_E))$ where $\rho$ denotes the spectral radius, i.e., the modulus of the eigenvalue of biggest modulus of the incidence matrix $A_E$ of the graph (see also \cite{MBMGRC2026}). In particular, this shows that the algebraic entropy of a Leavitt path algebra and of the associated ordinary path algebra coincide.

In this paper, we consider a closely related invariant, the growth series, defined as the generating function of the series $\dim_K\!\left( L_K(E)_n/L_K(E)_{n-1} \right)$, i.e.
$$
		\gr_{E}(t) = \sum_{n=0}^\infty \dim_K\!\left( L_K(E)_n/L_K(E)_{n-1} \right) t^n,
	$$
	where $L_K(E)_{-1} := \{0\}$. We show

\begin{teoInt}\label{thm:growth}
	Let $E$ be a finite graph with adjacency matrix $A_E$. Then the growth series of
	the Leavitt path algebra $L_K(E)$ is a rational function of the form
	$$\gr_E(t)=\frac{p(t)}{q(t)}\,\text{ with }\, p(t),q(t)\in\mathbb{Q}[t]$$
	so that in reduced form, $q(t)$ divides the square of the  reciprocal of the minimal polynomial of $A_E$.
    If $\rho(A_E)>0$, the power-series expansion of $\gr_E(t)$ converges for $|t|<1/\rho(A_E)$. In other case, $A_E$ is nilpotent and
	$\gr_E(t)$ is a polynomial.
\end{teoInt}

In particular, we see that every pole of the growth function is the reciprocal of a nonzero eigenvalue of the adjacency matrix $A_E$. Theorem \ref{thm:growth} also implies that the coefficients of $\gr_E(t)$ satisfy a linear recurrence and grow at most exponentially. This was already known, in fact the result mentioned above that $\entr(L_K(E))=\log(\rho(A_E))$ implies that the growth rate is determined by $\rho(A_E)$. In particular, we see that there are three possibilities: either the filtration $(L_K(E))_n$ stabilizes which happens precisely if the dimension of $L_K(E)$ is finite, or the growth of the coefficients of $\gr_E(t)$ is polynomial, or this growth is exponential. It is easy to characterize in terms of either the graph or equivalently of the adjacency matrix which is the case, we do that in Corollaries \ref{cor:growth-zero} and  \ref{cor:growth-polynomial}.


Given a  group $G$ with a finite generating system $X$, one can filter the group in terms of word length in $X$ and the size of the $n$-th level of this series yields the usual notion of group growth. For example, the growth of  the groups $F$, $T$ and $V$ that were introduced by Richard Thompson  in the $1960$s  (see Section \ref{sec:preThompson}) is known to be exponential, but the precise rate  is not known. However, in \cite{Burillo1}, Burillo is able to compute the precise growth rate for the monoid $F^+$ of positive words in $F$.

In this paper we consider a different notion of growth for graph topological full groups $\TFG$. These groups can be represented inside the Steinberg algebra and therefore also inside the Leavitt path algebra. In the case of Thompson's group $V$, which is the topological full group of the 2-rose, this fact goes back to Pardo \cite{Pardo}. The filtration of $L_K(E)$ can then be used to filter the topological full group, but it should be stressed that this is not a filtration by subgroups, nor the usual group filtration by word length above. In fact, it turns out  that in most cases this growth is much faster than the usual group growth. In Theorem \ref{theorem:TFGbound} we give a general upper bound on $|(L_K(E))_n\cap \TFG|$
and make it more precise in some examples.

In section \ref{sec:preThompson} we focus on the case of Thompson's group $V$. Computing the exact value of $|(\L_2)_n\cap V|$ where  $L_2$ is the Leavitt path algebra of the 2-rose is extremely difficult so we start by looking at the monoid $F^+$ consisting of elements that can be expressed as  positive words in the standard infinite generating system of $F$.



In Proposition \ref{prop:recurrencia} we give a recursive formula for the value of $|(\L_2)_n\cap F^+|$ and later we consider a finitely generated submonoid   $M$ 
 of $F^+$ to show in Theorem \ref{thm:Fplus-filtered-count} that for $M$ it is indeed possible to give a closed formula for $|(\L_2)_n\cap M|$ and that the growth of this series is exponential. It is worth to remark that the monoid $F^+$ was already considered in Thompson's handwritten notes, dated around 1973 or later in \cite{ThompsonNotes}. The submonoid $M$ appears in \cite[Theorem~4.6]{CannonFloydParry1996} where it is shown to be freely generated by $x_0$ and $x_1$.

 Things change dramatically when we consider $|(\L_2)_n\cap V|$. To bound this value, we introduce a series $C_k^n$ of combinatorial numbers which are related to the Catalan numbers. We give a recursive formula to compute them, alternatively, we show in Proposition \ref{prop:Ckn-generating-polynomial}
	that the generating functions 
	$$F_n(t):=\sum_{k\geq0}C_k^n t^k =\sum_{k=0}^{2^n-1}C_k^n t^k$$ satisfy
	
$$F_0(t)=1,\qquad F_{n+1}(t)=1+tF_n(t)^2,$$
a fact that can be used to compute small values of the $C_k^n$ numbers. Moreover, setting $x_n:=\sum_{k\geq0}C_k^n,$ we see that the series $x_n$ satisfies the  recurrence $x_{n+1}=1+x_n^2$. In fact, this series  has been considered in \cite{Aho} where the authors give a closed formula in terms of certain constant $\delta=1.502837$. Using this closed formula we prove

\begin{teoInt}\label{teoInt:doubleexp}
    There exists a constant $\delta=1.502837...$ such that for every $n\geq0$,
	\begin{equation}\label{eq:V-double-exponential-bounds}
		\left\lfloor\delta^{2^{\lfloor n/2\rfloor}}\right\rfloor
		\leq \left|V\cap(\mathbb L_2)_n\right|
		\leq \left\lfloor\delta^{3\cdot 2^n}\right\rfloor.
	\end{equation}
	\end{teoInt}

Moreover, we show that the lower bound also holds true for the group $F$, a fact that contrasts with our results for the monoid $M$ above: $F$ and $M$ are generated by the same two elements, but one as a group and the other as a monoid, in the first case the growth is double exponential, in the second only exponential.


%


The paper is organized as follows. We begin Section \ref{sec:Leavitt} by recalling the definition of Leavitt path algebra. We also introduce the standard filtration and some previous results on its growth and proceed to the proof of  Theorem \ref{thm:growth} and Corollaries \ref{cor:growth-zero} and  \ref{cor:growth-polynomial} where we characterize the cases when the filtration stabilizes or grows polynomially. In Section \ref{sec:TFG} we analyse the filtration induced in a graph topological full group by the standard filtration of the Leavitt path algebra. Again, we begin this section with some definitions, to introduce the basic notions and also to fix notation. We give a general upper bound for $|(L_K(E))_n\cap \TFG|$ and illustrate it with some examples. 
Section~\ref{sec:preThompson} is devoted to the particular case of Thompson's groups $V$ and $F$ and the monoids $F^+$ and $M$ described above. We prove some formulas of combinatorial nature and we finish by proving Theorem \ref{teoInt:doubleexp} where we give double exponential lower and upper bounds for $\left|V\cap(\mathbb L_2)_n\right|$.

 \medskip

\noindent{{\bf AI Disclosure:} No AI tool was used in the preparation of this work.}

\section{The growth series of Leavitt path \texorpdfstring{$K$}{K}-algebras of finite graphs}\label{sec:Leavitt}

\subsection{Preliminaries about Leavitt path \texorpdfstring{$K$}{K}-algebras}\label{sec:preliminaries}
This section provides the fundamental background about Leavitt path $K$-algebras. Here, $K$ is a field that we will assume fixed throughout the paper.

A \emph{directed graph} $E=(E^{0}, E^{1}, r, s)$ consists of two sets $E^{0}$ and $E^{1}$ whose elements are called \emph{vertices}
and \emph{edges} respectively and two functions $r, s : E^{1} \rightarrow E^{0}$ such that given an edge $e$ from
a vertex $v$ to a vertex $w$, $s(e)=v$ and $r(e)=w$. A vertex $v$ for which $s^{-1}(v)=\emptyset$ is called a sink, while a vertex $v$ for which $r^{-1}(v)=\emptyset$ is called a source. We will denote the set of sinks of $E$ by $\operatorname{Sink}(E)$ and the set of sources by $\operatorname{Source}(E)$. We say that a vertex $v\in E^0$ is an infinite emitter if $|s^{-1}(v)|=\infty$. The set of regular vertices (those which are neither sinks nor infinite emitters) is denoted by $\operatorname{Reg}(E)$. We will assume that $E$ is finite so throughout, graph will always mean finite directed graph.

A path of positive length is a word $\alpha=e_1\cdots e_m$ of edges satisfying $r(e_i)=s(e_{i+1})$ for $1\leq
	i<m$. Its length is $l(\alpha)=m$, its source is $s(\alpha)=s(e_1)$, and its range is $r(\alpha)=r(e_m)$. Vertices
are regarded as paths of length zero.

Graphs encode the combinatorial data from which the following algebra is constructed. Vertices and edges become algebraic generators and there is an extra type of generators that we call
\emph{ghost-edges} that can be understood as the reverses of the edges in $E^1$. 

\begin{definition}\label{def:LPA}\cite[Definition 1.2.3]{AAS}
	Let $E$ be an arbitrary (directed) graph and $K$ any field. We define a set $\left(E^{1}\right)^{*}$ consisting of symbols of the form $\left\{e^{*} \mid e \in\right.$ $\left.E^{1}\right\}$. The \emph{Leavitt path algebra} of $E$ with coefficients in $K$, denoted $L_{K}(E)$, is the free associative $K$-algebra generated by the set $E^{0} \cup E^{1} \cup\left(E^{1}\right)^{*}$, subject to the following relations:
	$$
		\begin{array}{ll}
			(V)   & v v^{\prime} = \delta_{v, v^{\prime}} v \quad \text{for all } v, v^{\prime} \in E^{0},                 \\
			(E1)  & s(e) e = e r(e) = e \quad \text{for all } e \in E^{1},                                                 \\
			(E2)  & r(e) e^{*} = e^{*} s(e) = e^{*} \quad \text{for all } e \in E^{1},                                     \\
			(CK1) & e^{*} e^{\prime} = \delta_{e, e^{\prime}} r(e) \quad \text{for all } e, e^{\prime} \in E^{1},          \\
			(CK2) & v = \sum_{\{ e \in E^{1} \mid s(e) = v \}} e e^{*} \quad \text{for every regular vertex } v \in E^{0}.
		\end{array}
	$$
    where $\delta_{x, x^{\prime}}$ is the Kronecker delta.
	\end{definition}


The algebra defined by the presentation with generators $E^{0}$, $E^{1}$ and only relators (V) and (E1) is called the \emph{path algebra} on $E$ and denoted by $KE$, while the algebra defined by the presentation with generators $E^{0}$, $E^{1}$, $\left(E^{1}\right)^{*}$ and relatiors (V), (E1) and (E2) is called the \emph{path algebra of the extended graph}  $\Hat{E}$ and denoted by $K\Hat{E}$. 

An important example that motivated the theory of Leavitt path algebras is the case when the graph is the $r$-rose $R_r$, that is, one vertex and $r$ loops. This algebra is known as the
Leavitt algebra of type $(1,r)$ and was introduced by Leavitt in connection with the study
of the invariant basis number (IBN) property.

\begin{definition}\label{def:typen}
	Let $r\geq2$, the \emph{Leavitt algebra of type $(1,r)$} is
	$$
		\mathbb L_r=
		K\langle y_1,\ldots,y_r,y_1^*,\ldots,y_r^*
		\mid y_i^*y_j=\delta_{i,j},\
		\sum_{i=1}^r y_iy_i^*=1\rangle .
	$$
    Sending the $i$-th loop of $R_r$ to $y_i$ and its ghost edge to $y_i^*$ yields an isomorphism $L_K(R_r)\cong\mathbb L_r$.
\end{definition}

\subsection{The standard filtration.}
Relators (V), (E1) and (E2) in Definition \ref{def:LPA} imply that elements of a Leavitt path algebra are linear combinations of paths. This means that following \cite[Definition 3.10]{entropy1}, path lengths can be used to define a natural filtration of the algebra and an associated notion of growth.  

\begin{definition}\label{def:filtLPA}
	Let $L_{K}(E)$ be the Leavitt path $K$-algebra of the graph $E$. For each $n\geq 0$ let $L_{K}(E)_{n}$ be the linear span of the set of monomials $\lambda\mu^{*}$ with $l(\lambda)+l(\mu)\leq n$ (recall that vertices have length zero). Then $\{L_{K}(E)_{n}\}_{n\in \mathbb{N}}$ is the \textit{standard filtration} of $L_{K}(E)$.\end{definition}

With respect to this filtration, $L_K(E)$ is a \emph{filtered algebra}. Recall that a filtered algebra is a $K$-algebra $A$ with a series of subspaces
$$A_1\subseteq A_2\subseteq\ldots\subseteq A_n\subseteq\ldots$$
such that $A_iA_j\subseteq A_{i+j}$ and $A=\cup A_n$. We are going to consider a special type of filtered algebras where $A_1$ is a finite dimensional generating set and $A_i$ is the $K$-span of the products $a_1\ldots a_i$ for $a_1,\ldots a_i\in A_1$. This implies that if $A_1$ is finite dimensional, then so is each $A_i$. A way to measure the growth of the series is using the \emph{Gelfand-Kirillov} dimension, defined as
$$\GK(A)=\limsup_{n\to\infty}\frac{\log(\dim_KA_n)}{\log(n)}.$$ 
This dimension does not depend on the chosen generating system $A_1$ as long as it is finite dimensional. Moreover, the Gelfand-Kirillov dimension is finite if and only if the filtration growth is polynomially bounded (\cite{Zelmanov}).
In \cite{Zelmanov}, this condition is characterized for the case in which $A$ is a Leavitt path algebra:

\begin{theorem}\cite[Theorem 5]{Zelmanov}\label{thm:Zelmanov} Let $E$ be a finite graph. Then $\GK(L_K(E))$ is finite if and only if $E$ has the property that any two distinct cycles of $E$ do not have a common vertex.
\end{theorem}

The condition on the graph in the previous theorem is rather restrictive, so one may ask if there is a finer way to measure growth for arbitrary Leavitt path algebras. To do that, 
in \cite{entropy1} the authors introduce the notion of  \emph{algebraic entropy} of a filtered algebra $A$ with respect to the filtration $\mathcal{F}=(A_n)$ as $\entr(A,\mathcal{F})=0$ if $A$ has finite dimension or  
$$\entr(A,\mathcal{F})=\limsup_{n\to\infty}\frac{\log(\dim_K(A_n/A_{n-1}))}{n}$$ 
otherwise. The algebraic entropy does depend on the precise filtration, but when we consider Leavitt path algebras we will fix the filtration of definition \ref{def:filtLPA} and simplify notation to $\entr(L_K(E))$.

In the main result of \cite{entropy2}, the authors determine $\entr(L_K(E))$ in terms of the \emph{adjacency matrix} of the graph, defined as follows. Let $E$ be a finite graph with $d$ vertices which are labeled from 1 to $d$. The adjacency matrix $A_E$ is the  $d\times d$-matrix having as $(i,j)$-entry the number of edges in $E^1$ with source the vertex $i$ and range the vertex $j$. For any matrix $A$ we denote by $\rho(A)$ the spectral radius, i.e., the maximum of the modulus of the eigenvalues of $A$. Bock, Gil Canto, Mart\'in Barquero,
Mart\'in Gonz\'alez, Ruiz Campos and Sebandal show:

\begin{theorem}\cite[Theorem 3.5]{entropy2}\label{thm:entropy} Let $E$ be a finite graph. Then 
$$
\entr(L_K(E))=\begin{cases} 0 &\text{if $E$ is acyclic},\\
\log(\rho(A_E)) &\text{otherwise}.
\end{cases}
$$
\end{theorem}


A key part of the argument is the following computation which we state as in \cite{entropy1} but comes from a previous paper by Bock and Sebandal \cite[Theorem~38]{BockSebandal}. 
It is based on the fact that given two vertices $i,j$ and some $\alpha\geq 0$, the number of paths of length $\alpha$ with source $i$ and range $j$ is precisely $(A^\alpha_E)_{i,j}$.

\begin{prop}\label{prop:tecentropy}\cite[Proposition 6.1]{entropy1}
	Let $E$ be a finite graph with $d$ vertices, with adjacency matrix $A_E$ and fix $n\in\mathbb{N}$. Then,
	\begin{equation}
		\dim\!\left(L_K(E)_n/L_K(E)_{n-1}\right)
		=
		\sum_{s=0}^{n}\sum_{j=1}^{d}
		\left(\sum_{i=1}^{d}(A_E^s)_{i,j}\right)
		\left(\sum_{i=1}^{d}(A_E^{n-s})_{i,j}\right)
		-\sum_{s=1}^{n-1}\sum_{j=1}^{d}
		\varepsilon_j
		\left(\sum_{i=1}^{d}(A_E^{s-1})_{i,j}\right)
		\left(\sum_{i=1}^{d}(A_E^{n-s-1})_{i,j}\right),
	\end{equation}
	where
	$$
		\varepsilon_j=
		\begin{cases}
			0, & \text{ if $j$ is a sink} \\[6pt]
			1, & \text{otherwise}.
		\end{cases}
	$$
\end{prop}
Here, we say that a vertex $j$ is a \emph{sink} if it is not the source of any arrow. Equivalently, if
$\sum_{m=1}^{d}(A_E)_{j,m}=0.$

\subsection{The growth series.}\label{sec:growth-series}

 In this section we prove Theorem \ref{thm:growth}.

  Given a finite graph $E$, we define the  \emph{growth series} of the Leavitt path algebra $L_K(E)$  associated to the filtration of Definition \ref{def:filtLPA}  as the generating function 
	$$
		\gr_{E}(t) = \sum_{n=0}^\infty \dim_K\!\left( L_K(E)_n/L_K(E)_{n-1} \right) t^n,
	$$
	where $L_K(E)_{-1} := \{0\}$. We also set

    $$\GK_E(t)= \sum_{n=0}^\infty \dim_K\! L_K(E)_n t^n.$$

Observe that these two series are related as follows
\begin{equation}\label{eq:grGK}
\begin{aligned}
    \gr_{E}(t) &= \sum_{n=0}^\infty (\dim_K\!L_K(E)_n-\dim_K\!L_K(E)_{n-1}) t^n\\
 &=\sum_{n=0}^\infty \dim_K\!L_K(E)_n t^n- \sum_{n=0}^\infty \dim_K\!L_K(E)_{n-1} t^n\\
 &=\GK_E(t)-t\GK_E(t)=(1-t)\GK_E(t).
\end{aligned}
\end{equation}

\bigskip

 
 \noindent{\sl Proof of Theorem \ref{thm:growth}.} Recall that given two vertices $i,j$, $(A_E^\alpha)_{ij}$ counts the number of paths of length $\alpha$ from $i$ to $j$. Therefore the number of paths of length $\alpha$ which end at $j$ is
	\begin{equation}\label{eq:tau}
		\tau_j(\alpha):=\sum_{i=1}^{d}(A_E^\alpha)_{ij}.
	\end{equation}

	As in Proposition \ref{prop:tecentropy}, let 
	$\varepsilon_j$ be 0 if $v_j$ is a sink and 1 otherwise and set
	\begin{equation}\label{eq:Tau}
		T_j(t):=\sum_{\alpha=0}^{\infty}\tau_j(\alpha)t^\alpha.
	\end{equation}

	If we write the expression of Proposition~\ref{prop:tecentropy} using \eqref{eq:tau} and \eqref{eq:Tau}, the
	definition of the growth series gives
	\begin{equation}\label{eq:claim1-path-count}
		\begin{aligned}
			\gr_E(t)
			 & =\sum_{n=0}^{\infty}\left(
			\sum_{s=0}^{n}\sum_{j=1}^{d}\tau_j(s)\tau_j(n-s)
			-\sum_{s=1}^{n-1}\sum_{j=1}^{d}
			\varepsilon_j\tau_j(s-1)\tau_j(n-s-1)
			\right)t^n                                      \\
			 & =\sum_{j=1}^{d}T_j(t)^2
			-t^2\sum_{j=1}^{d}\varepsilon_jT_j(t)^2          =\sum_{j=1}^{d}(1-\varepsilon_jt^2)T_j(t)^2,
		\end{aligned}
	\end{equation}
	where the last equality is obtained as follows: the first term is a combination of the Cauchy product formula and \eqref{eq:Tau}, and for the second,
	put $a=s-1$ and $b=n-s-1$. Then $a,b\geq0$ and $n=a+b+2$, which gives
	$t^2T_j(t)^2$ after interchanging the order of summation. At this point, to deduce that the functions $T_j(t)$ are rational functions with denominator in reduced form dividing the reciprocal of the characteristic polynomial of $A_E$ we can use well known arguments, see for example \cite[Proposition V.6]{FS}, but we include here a proof explicitly that yields the reciprocal of the minimal polynomial instead.

	Let $\mathbf{1}$ be the column vector all of whose entries are $1$, and let $e_j$ be the $j$-th vector in the canonical
	basis of $K^d$. Then
	$$
		\tau_j(\alpha)=\mathbf{1}^{\mathsf T}A_E^\alpha e_j.
	$$
	As a formal matrix series, we obtain
	$$
		\begin{aligned}
			T_j(t)  =\sum_{\alpha=0}^{\infty}\tau_j(\alpha)t^\alpha                          =\sum_{\alpha=0}^{\infty}\mathbf{1}^{\mathsf T}A_E^\alpha e_j\,t^\alpha =\mathbf{1}^{\mathsf T}(I-tA_E)^{-1}e_j,
		\end{aligned}
	$$
	by the Neumann series identity in $\mathbb{Q}[[t]]$. If $\rho(A_E)>0$, the corresponding analytic series converges
	for $|t|<1/\rho(A_E)$. In other case, $A_E$ is nilpotent, the matrix series is finite and there is no analytic restriction.

	Now, we express $A_{E}$ in its Jordan form. Thus there is an invertible matrix $P$ with entries in $\mathbb{C}$
	such that $A_E=PJP^{-1}$, with $J$ in Jordan form. Hence
	$$
		(I-tA_E)^{-1}=P(I-tJ)^{-1}P^{-1}.
	$$
	A Jordan block of size $m$ has the form $J_{\lambda,m}=\lambda I+N$, where $N^m=0$. Since
	$$
		I-tJ_{\lambda,m}=(1-\lambda t) \left(I-\frac{t}{1-\lambda t}N\right),
	$$
	and the geometric series in the nilpotent matrix $N$ stops after its $(m-1)$-st power, its inverse can be written
	as
	\begin{equation}\label{eq:diagblocks}
		(I-tJ_{\lambda,m})^{-1}=\frac{1}{1-\lambda t}(I-\frac{t}{1-\lambda t}N)^{-1}
		=\sum_{r=0}^{m-1}\frac{t^rN^r}{(1-\lambda t)^{r+1}}.
	\end{equation}
        
         Let $\Lambda$ be the set of distinct nonzero eigenvalues of $A_E$, and for $\lambda\in\Lambda$ let $m_\lambda$ be
	the size of the biggest Jordan block associated with $\lambda$. Let also
	$$
		D(t):=\prod_{\lambda\in\Lambda}(1-\lambda t)^{m_\lambda}
	$$
    be the reciprocal of the minimal polynomial of the matrix $A_E$. The blocks in (\ref{eq:diagblocks}) are the diagonal blocks of $(I-tJ)^{-1}$ and multiplication on the left by $\mathbf{1}^{\mathsf T}P$ and on the right by $P^{-1}e_j$ only replaces
	the matrices $N^r$ by scalar coefficients. This means that combining the contributions of all the Jordan blocks, one gets
	$$	T_j(t)=\frac{h_j(t)}{D(t)}$$
    for certain polynomial $h_j(t)\in\mathbb{C}[t]$. 
	Substitution in \eqref{eq:claim1-path-count} implies that $\gr_E(t)$ is rational and that in reduced form the denominator divides $D(t)^2$.
    Since
	$$
		(I-tA_E)^{-1}
		=\frac{\operatorname{adj}(I-tA_E)}{\det(I-tA_E)},
	$$
	and $A_E$ has integer entries, we see that the polynomials $N_j(t)$ and $D(t)$ have coefficients in $\Q$. 
    Finally, the  assertion about convergence follows from the Neumann series and the observation that $\rho(A_E)=0$ if and only if
	 $A_E$ is nilpotent, and in that case $\gr_E(t)$ is a polynomial.

\begin{remark}\label{rem:AnCom}
    In general, see for example \cite[Theorem IV.7]{FS}, if $f(t)$ is a complex function analytic at 0, the radius $R$ of convergence, i.e., the smallest modulus of a singularity, measures the growth of the coefficients of the Taylor expansion of $f$, in the sense that this growth is equivalent to the growth on $n$ of $1/R^n$. This means that for a rational function, the growth is exponential on the inverses of the poles of the function. In our case, these poles are the inverses of the eigenvalues of the adjacency matrix. In fact, the particular case that we are interested in follows from the easy fact that
    $$\left(\frac{1}{1-t\lambda}\right)^m=\sum_{n=0}^\infty a_n\lambda^nt^n$$
    where $a_n=\frac{(n+m-1)!}{n!(m-1)!}$ is a polynomial in $n$. So if $\|\lambda\|>1$ we get exponential growth and if $\|\lambda\|=1$, the growth is polynomial.
\end{remark}

We illustrate with two examples how one can compute more precisely the growth polynomial $gr_E(t)$.  

\begin{example}\label{ex:Vgrowth} Let $R_2$ be the 2-rose, i.e. the graph
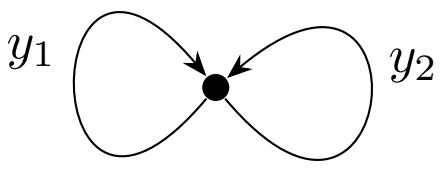
\begin{figure}[H]
\centering
\begin{tikzpicture}[
    vertex/.style={fill,circle,inner sep=1.8pt},
    >={Stealth}
]
    \node[vertex] (v) at (0,0) {};

    \draw[->] (v)
        to[loop,looseness=12,min distance=18mm,in=130,out=-130]
        node[midway,above left] {$y_1$} (v);

    \draw[->] (v)
        to[loop,looseness=12,min distance=18mm,in=40,out=-50]
        node[midway,above right] {$y_2$} (v);
\end{tikzpicture}
\end{figure}
Its Leavitt path algebra  is the Leavitt algebra of type $(1,2)$ and has presentation
$$
\mathbb L_2:=L_K(R_2)= K\!\left\langle y_1,y_2,y_1^*,y_2^*\ \middle|\
y_i^*y_j=\delta_{ij},\quad y_1y_1^*+y_2y_2^*=1,
\right\rangle.
$$
Using~\eqref{eq:claim1-path-count} in the proof of Theorem~\ref{thm:growth} we have
$$
\tau_1(\alpha)=2^{\alpha},\quad T_1(t)=\sum_{\alpha=0}^{\infty}(2t)^{\alpha}=\frac{1}{1-2t},\,\,|t|<\frac{1}{2}
$$
and thus
$$
gr_{R_2}(t)=(1-t^2)\left(\frac{1}{1-2t}\right)^2.
$$
This  series converges analytically for $|t|<\frac{1}{2}$, with $2=\rho(A_E)$ (Theorem \ref{thm:growth}).
Alternatively, we could use the explicit formula in \cite[Example 3.12]{entropy1} for $\dim_K\left( L_K(R_2)_n/L_K(R_2)_{n-1} \right)=2^{n-2}(3n+5)$, or compute 
$GK_{R_2}(t)$ and use ~\eqref{eq:grGK}. 
\end{example}


\begin{example}\label{ex:Toeplitzgrowth} Let $E_T$ be the graph
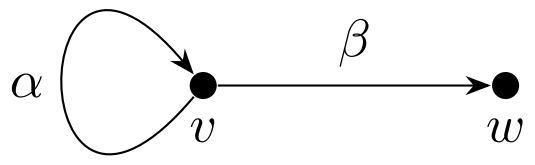
\begin{figure}[H]
\centering
\begin{tikzpicture}[vertex/.style={fill,circle,inner sep=1.8pt},>={Stealth}]
\node[vertex,label=below:{$v$}] (v) at (0,0) {};
\node[vertex,label=below:{$w$}] (w) at (2,0) {};
\draw[->] (v)
    to[loop,looseness=12,min distance=18mm,in=130,out=-130]
    node[midway,left] {$\alpha$} (v);
\draw[->] (v) --
    node[midway,above] {$\beta$} (w);
\end{tikzpicture}
\end{figure}
Its Leavitt path algebra is the Toeplitz algebra
$$
L_K(E_T)\cong K\langle x,y\mid xy-1\rangle =K\langle x,y\rangle/\langle xy-1\rangle.
$$
The adjacency matrix  is 
$$
A_{E_T}=\begin{pmatrix}
1&1\\
0&0 
\end{pmatrix}=A_{E_T}^{\alpha},\,\, \alpha\geq 1.
$$
Using again~\eqref{eq:claim1-path-count}, $\tau_1(\alpha)=\tau_2(\alpha)=1$ for all $\alpha\geq 0$ and hence $T_1(t)=T_2(t)=\frac{1}{1-t}$ for $|t|<1$ thus
$$
gr_{E_T}(t)=(1-t^2)\left(\frac{1}{1-t}\right)^2+\left(\frac{1}{1-t}\right)^2=\frac{2-t^2}{(1-t)^2}.
$$

\end{example}

We can also deduce the following consequence.

\begin{cor}\label{cor:growth-zero}
	Let $E$ be a finite graph, $A_{E}$ its adjacency matrix and $\gr_{E}(t)$ the growth series of its associated Leavitt path $K$-algebra. Then, the following are equivalent:
	\begin{enumerate}
		\item[i)] $\gr_{E}(t)\in \mathbb{Z}_{\geq 0}[t]$,
		\item[ii)] $L_{K}(E)$ is a finite-dimensional Leavitt path $K$-algebra,
		\item[iii)] $E$ is acyclic,
		\item[iv)] $A_{E}$ is nilpotent,
	\end{enumerate}

\end{cor}
\begin{proof}
	For each $n\geq 0$, the coefficient
	$$
		h_n:=\dim_K\!\left(L_K(E)_n/L_K(E)_{n-1}\right)
	$$
	is a nonnegative integer. Since the standard filtration is increasing, $\gr_E(t)$ is a polynomial if and only if
	$h_n=0$ for all sufficiently large $n$, that is, if and only if the filtration stabilizes. This is equivalent to
	$L_K(E)$ being finite-dimensional. It is clear that in that case $\gr_E(t)\in\mathbb{Z}_{\geq 0}[t]$.

	By \cite[Theorem 2.6.17]{AAS}, $L_K(E)$ is finite-dimensional if and only if $E$ is acyclic.  Now, let $d=|E^0|$ and let $E$ be acyclic, then there are no paths of
	length $d$, and thus $A_E^d=0$. Conversely, if $E$ contains a cycle of length $r$ based at $v_i$, then
	$(A_E^{rt})_{ii}>0$ for any $t\geq 1$, and hence $A_E$ is not nilpotent.

 \end{proof}

\begin{cor}\label{cor:growth-polynomial}
	Let $E$ be a finite graph, $A_{E}$ its adjacency matrix and $\gr_{E}(t)$ the growth series of its associated Leavitt path $K$-algebra. Then, the following are equivalent:
	\begin{enumerate}
		\item[i)]  there are some positive integers $k,l$ such that $(1-t^k)^l\gr_E(t)\in \mathbb{Z}[t]$,
		\item[ii)] the series $\dim_K(L_{K}(E))_n$ grows polynomially,
		\item[iii)] $E$ has the property that any two distinct cycles do not have a common vertex,
		\item[iv)] the non-zero eigenvalues of $A_{E}$ are roots of unity.
	\end{enumerate}

\end{cor}
\begin{proof}

As $\gr_E(t)=(1-t)\GK_E(t)$, item  i) is equivalent to the same statement but with the function $\GK_E(t)$ and using Remark \ref{rem:AnCom}
this is equivalent to ii) (and also equivalent to the coefficients of $\gr_E(t)$ growing polynomially).
The equivalence between ii) and iii) is Theorem \ref{thm:Zelmanov}. 

The fact that iv) implies i) follows from Theorem \ref{thm:growth}. For the converse, we have already established that i) is equivalent to the fact that the coefficients of $\gr_E(t)$ grow polynomially. Taking into account Theorem \ref{thm:entropy}, this implies that $\rho(A_E)=1$ and therefore all the eigenvalues of $A_E$ have modulus at most one. As they are roots of the characteristic polynomial of $A_E$, which has integer coefficients, this implies by a well known result by Kronecker \cite{Kronecker} that they must all be either zero or roots of unity so we get iv).

  \end{proof}

\begin{remark} Instead of using Theorem \ref{thm:entropy}, one can check that iii) implies iv) directly as follows. Condition iii) implies that we may group the vertices of $E^0$ in disjoint sets $C_1,\ldots,C_k$ so that each $C_i$ consists of either a single vertex not lying in any cycle or the vertices of a cycle. If we chose the labels so that for $i\neq j$, if there is a path from a vertex in $C_i$ to a vertex in $C_j$, then $i<j$ and form the matrix $A_E$ according to the induced ordering of the vertices, then $A_E$ is an upper triangular block matrix such that the diagonal block $D_i$ is either zero if it corresponds to $C_i$ with a single vertex or satisfies $D_i^r=I$ if it corresponds to a cycle $C_i$ of length $r$. The eigenvalues of such a matrix can only be either zero or a root of 1.  
\end{remark}

\section{Graph topological full groups and the filtration induced from the Leavitt path algebra embedding.}\label{sec:TFG}

 As in Section \ref{sec:preliminaries}, let $E=(E^0,E^1,r,s)$ be a finite directed graph. Let $E^f$ and $E^\infty$ be the sets of finite and infinite paths in $E$ respectively and put  
$$\partial E=E^\infty\cup\{\alpha\in E^f\mid r(\alpha)\text{ is a sink}\}.$$

\subsection{The Renault-Deaconu groupoid of a graph and the associated topological full group}

The \emph{Renault-Deaconu groupoid}, or just \emph{graph groupoid} associated to $E$ is
$$\mathcal{G}_E=\{(\alpha x,l(\alpha)-l(\beta),\beta x)\mid \alpha,\beta\in E^f,x\in\partial E, r(\alpha)=r(\beta)=s(x)\}.$$
The product structure is given by $(\alpha,n,\beta)(\beta,m,\gamma)=(\alpha,n+m,\gamma)$, the set of units $\mathcal{G}_E^{(0)}$ consists of elements of the form $(\alpha,0,\alpha)$ for $\alpha\in\partial E$ and the inverse of $(\alpha,n,\beta)$ is $(\beta,-n,\alpha)$.
In \cite{Matui} and \cite{Nek} Matui and Nekrashevich give related but slightly different definitions for the topological full group of a groupoid. Here, we will follow Nekrashevich approach and consider only the case of graph groupoids as $\mathcal{G}_E$.

We define maps $s,r:\mathcal{G}_E\to\partial E$ given by $s(\alpha,n,\beta)=\alpha$, $r(\alpha,n,\beta)=\beta$. 
A \emph{bisection} in $\mathcal{G}_E$ is a subset $U\subseteq\mathcal{G}_E$ such that the maps $s,r$ are injective, if they are also surjective, then we say that the bisection is \emph{full}. The  groupoid $\mathcal{G}_E$ has a topology with a basis consisting of compact open bisections, which we describe next,  using a notation close to Section 9 in \cite{NylandOrtega}. Given $\mu,\nu\in E^f$ with $r(\mu)=r(\nu)$, let
$$Z(\mu,\nu)=\{(\mu x,l(\mu)-l(\nu),\nu x)\in\mathcal{G}_E \mid \mu,\nu\in E^{f}, x\in\partial E, r(\mu)=r(\nu)=s(x)\}.$$
For $\mu,\nu$ as before, let $F\subseteq r(\mu)E^1$. In \cite{NylandOrtega}, this set $F$ is required to be finite, but with our finiteness hypothesis we do not need to ask that. Now, set
$$Z(\mu,F,\nu)=Z(\mu,\nu)\setminus\cup_{e\in F}Z(\mu e,\nu e).$$
In $\mathcal{G}_E$, we consider the topology generated by these basic sets. It turns out that then these basic sets are compact open bisections and from this, one can deduce that $\mathcal{G}_E$ is an ample Hausdorff groupoid. In our case, as we are assuming that $E$ is finite, we have  
$$Z(\mu,F,\nu)=Z(\mu,\nu)\setminus\cup_{e\in F}Z(\mu e,\nu e).$$
Therefore, the compact open sets of the form $Z(\mu,\nu)$ are also a basis for the topology of $\mathcal{G}_E$. Now, if $U$ is a compact open full bisection, then it must be a union of those open sets and compactness implies that it is a finite union 
$$U=\cup_{i=1}^kZ(\mu_i,\nu_i).$$
Pairs of sets of the form $Z(\mu_i,\nu_i)$, $Z(\mu_j,\nu_j)$ are either disjoint or one is inside the other and this means that we may assume that the union above is a disjoint union. Then, the fact that $U$ is a full bisection implies that
$\partial E=s(U)$ is a disjoint union of \emph{cylinders}, i.e., of  sets of the following form
$$s(Z(\mu_i,\nu_i))=Z(\mu_i):=\{\mu_ix\mid x\in\partial E, r(\mu_i)=s(x)\}$$
and the same for $r(U)$ and the sets $Z(\nu_i)$. As this property will be important later, we will give such sets a name.

\begin{definition}
    Let $\mu_1,\ldots,\mu_k$ be a family of elements in $E^{f}$,   such that 
$$\partial E=\sqcup_{i=1}^kZ(\mu_i),$$
Then we say that $\{\mu_1,\ldots,\mu_k\}$ is a \emph{forest-basis} of $\partial E$.
\end{definition}
Observe that if $\{\mu_1,\ldots,\mu_k\}$ is a forest-basis of $\partial E$, then $\{s(\mu_1),\ldots,s(\mu_k)\}=E^0$. Given two compact open full bisections $U$, $V$, its product
$$\{gh\mid g\in U,h\in V\mid r(g)=s(h)\}$$
is also a compact open full bisection. The set of units can be seen as a trivial compact open full bisection:
$$\mathcal{G}_E^{(0)}=\{(x,0,x)\mid x\in\partial E\}$$
and given a compact open full bisection $U$, if we set
$$U^{-1}=\{(\alpha,n,\beta)\mid (\beta,-n,\alpha)\in U\},$$
then $UU^{-1}=U^{-1}U=\mathcal{G}_E^{(0)}$.
We define the \emph{topological full group} $\TFG$ of $E$ as the set of compact open full bisections 
$$\TFG:=\{U\subseteq\mathcal{G}_E\text{ compact open full bisection}\}$$
with the product above. Every $U\in\TFG$ induces a homeomorphism of $\partial E$ given by
$$
\pi_U=s|_U\circ(r|_U)^{-1}.
$$
Indeed, both restrictions are homeomorphisms onto $\partial E$. Moreover,
$$
\pi_{UV}=\pi_U\circ\pi_V,
\qquad
\pi_{\mathcal{G}_E^{(0)}}=\operatorname{id}_{\partial E},
$$
so $U\cdot x:=\pi_U(x)$ defines an action of $\TFG$ on $\partial E$. Passing from the bisections to the homeomorphisms they induce, leads to Matui's convention for the topological full group, denoted by $[[\mathcal{G}_{E}]]$; see \cite[Definition 2.1]{Matui}.
$$
[[\mathcal{G}_E]]
:=
\{\pi_U\mid U\in\TFG\}
\leq \operatorname{Homeo}(\partial E),
$$
with composition as its group operation.

\begin{remark}\label{rem:L}
There are two related conventions for the topological full group. As we have mentioned, Nekrashevych's $\TFG$ consists of the compact open full bisections, whereas Matui's $[[\mathcal{G}_E]]$ consists of the homeomorphisms of $\partial E$ induced by such bisections. So the action above gives a group epimorphism
$$
\rho:\TFG\longrightarrow[[\mathcal{G}_E]],
\qquad U\longmapsto\pi_U.
$$

This homomorphism is surjective, Hence
$$
[[\mathcal{G}_E]]\cong\TFG/\ker\rho.
$$
And $\rho$ is injective exactly when the actions of two different bisections are always different. This can be made precise using the \emph{isotropy} of the groupoid, defined as
$$
\operatorname{Iso}(\mathcal{G}_E)=\{g\in\mathcal{G}_E\mid s(g)=r(g)\}.
$$
A full bisection acts trivially on $\partial E$ when $U\subseteq\operatorname{Iso}(\mathcal{G}_E)$. As $U$ is open, this is equivalent to $U\subseteq\operatorname{Iso}(\mathcal{G}_E)^\circ$ where the superscript $\circ$ denotes the interior.  Therefore $\rho(U)=\rho(V)$ if and only if $V^{-1}U\subseteq\operatorname{Iso}(\mathcal{G}_E)^\circ$ thus
$$
\ker\rho=\{U\in\TFG\mid U\subseteq\operatorname{Iso} (\mathcal{G}_E)^\circ\}.
$$
When $\rho$ is injective, i.e., when $\operatorname{Iso}(\mathcal{G}_E)^\circ=\mathcal{G}_E^{(0)}$ the groupoid is called \emph{effective}; see \cite[Definition~2.2]{Nek}. This can be characterized in terms of the graph:   $\mathcal{G}_E$ is effective if and only if the graph has condition (L), meaning that every cycle has an exit (a recollection of this fact can be found in \cite[Section 3.3]{NylandOrtegaAH}).
\end{remark}

To describe the elements of the group $\TFG$, observe that by the discussion above, any $U\in\TFG$ can be written as a disjoint union
\begin{equation}\label{eq:bisection}U=\sqcup_{i=1}^kZ(\mu_i,\nu_i)\end{equation}
so that $B_1=\{\mu_1,\ldots,\mu_k\}$ and $B_2=\{\nu_1,\ldots,\nu_k\}$ are forest-bases. In fact, $U$ is determined by $B_1,B_2$ and a bijection $\varphi:B_1\to B_2$ so we will say that the triple $(B_1,B_2,\varphi)$ is a representative of $U$ by forest-bases. But the converse is not true, in the sense that, given $U$, the forest-basis $B_1$ and $B_2$ are not unique: to see it observe that if $\mu,\nu\in E^{f}$  have $r(\mu)=r(\nu)=v$ and $\{e_1,\ldots,e_t\}$ are the edges in $E^1$ having $v$ as source, then  
$$Z(\mu,\nu)=\sqcup_{i=1}^{k}Z(\mu e_i,\nu e_i).$$
Given $U$, an associated forest-bases representative for which the number $k$ is called a \emph{reduced} forest-basis representative of $U$.


\begin{example}\label{ex:V} It was first observed by Nekrashevych that Thompson's group $V$ is the topological full group of the 2-rose $R_2$, i.e., the graph with a single vertex and two loops that we have considered in Example \ref{ex:Vgrowth}.
 Then, we may represent forest-bases simply as finite rooted binary trees and group elements as bijections between them, see Section \ref{sec:preThompson} for details. Similarly, any Higman-Thompson group $V_{d,n}$ is the topological full group of the graph $E_{d,n}$ obtained by attaching an $n$-rose to the sink of a directed $n$-line. 
In this case reduced forest-bases are unique as we may identify forest-bases representatives with pairs of reduced binary trees (see Section \ref{sec:preThompson}).
\begin{figure}[H]\label{E32}
  \centering
    \begin{tikzpicture}[scale=0.75, transform shape,vertex/.style={fill,circle,inner sep=1.8pt},>={Stealth} ]
        \node[vertex] (v1) at (0,0) {};
        \node[vertex] (v2) at (2,0) {};
        \node[vertex] (v3) at (4,0) {};

        \draw[->] (v1) -- (v2);
        \draw[->] (v2) -- (v3);

        \draw[->] (v3)
            to[loop,looseness=12,min distance=18mm,in=130,out=-130] (v3);
        \draw[->] (v3)
            to[loop,looseness=12,min distance=18mm,in=40,out=-50] (v3);
    \end{tikzpicture}
    \caption{The graph $E_{3,2}$.}
\end{figure}
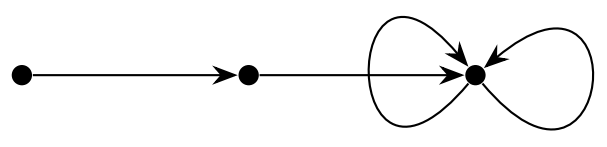
\end{example}

\begin{example} Consider now the graph  $E_{T}$ of example \ref{ex:Toeplitzgrowth}.
We want to describe the topological full group $\mathbf{F}(\mathcal{G}_{E_T})$. To do that, observe first that all the forest-bases come from iterated expansions of the vertex $v$ and those of its children having label also $v$. (see Figure~\ref{fig:forestbasis} below)
\begin{figure}[H]
\centering
\begin{tikzpicture}[
    black vertex/.style={
        circle,
        fill=black,
        draw=black,
        inner sep=1.1pt
    },
    white vertex/.style={
        circle,
        fill=white,
        draw=black,
        inner sep=1.1pt
    },
    level distance=8mm,
    sibling distance=7mm,
    >={Stealth}
]
\node[black vertex,label=below:{$v$}] (v) at (-3,0) {};
\node[white vertex,label=below:{$w$}] (w) at (-1.5,0) {};
\draw[->,red] (v)
    to[loop,looseness=12,min distance=18mm,in=130,out=-130]
    node[midway,left,text=red] {$\alpha$} (v);
\draw[->,blue] (v) --
    node[above,text=blue] {$\beta$} (w);
\node[black vertex] at (1.5,1) {}
    child {
        node[black vertex] {}
        child {
            node[black vertex] {}
            child {
                node [inner sep=1.1pt,label=below:{$\vdots$}] {}
                edge from parent[->,red,dashed]
            }
            child {
                node[white vertex] {}
                edge from parent[->,blue]
            }
            edge from parent[->,red]
        }
        child {
            node[white vertex] {}
            edge from parent[->,blue]
        }
        edge from parent[->,red]
    }
    child {
        node[white vertex] {}
        edge from parent[->,blue]
    };
\node[white vertex] at (3,1) {};
\end{tikzpicture}
\caption{The graph \(E_T\) and the trees of paths rooted at \(v\) and \(w\), respectively.}
\label{fig:forestbasis}
\end{figure}
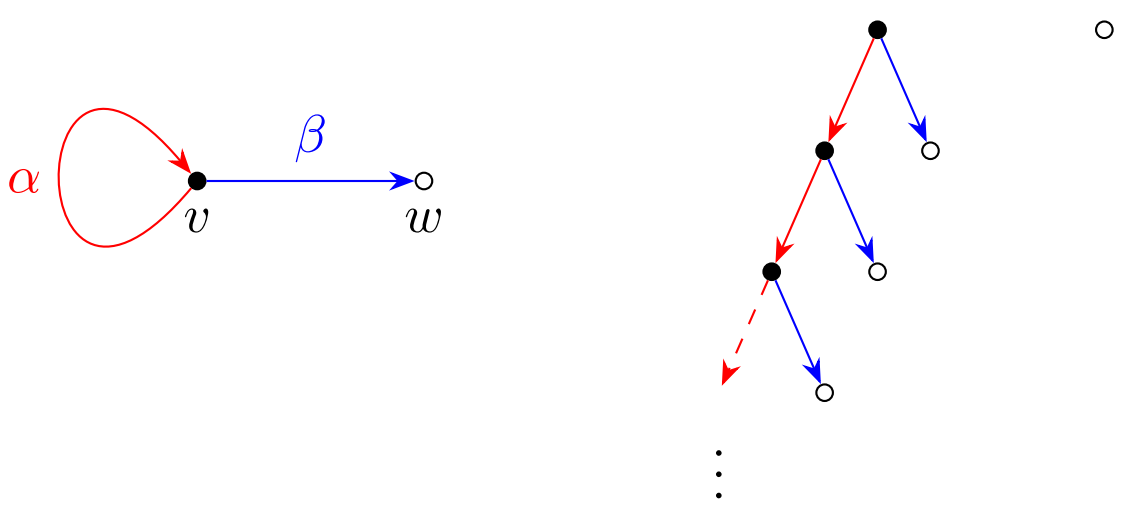
This means that there is only one forest-basis of each cardinality and all of them have the form
$$
B_n=\{\alpha^n,w\}\cup\{\alpha^j\beta:0\leq j\leq n-1\},\qquad n\geq0,
$$
where $\alpha^0=v$. So we have $|B_n|=n+2$ and in $B_n$ there is only one element, $\alpha^n$, with range $v$, and the other $n+1$ elements have range $w$.
Therefore the forest-bases pair expression of every $g\in\mathbf{F}(\mathcal{G}_{E_T})$ must be of the form $(B_n,B_n,\varphi)$ where $\varphi$ is a permutation between the $n+1$ elements having range $w$ and we get a copy of the symmetric group $\mathrm{Sym}(n+1)$, moreover expanding the forest-bases at $v$ see that these copies are nested into each other and we have
$$\mathbf{F}(\mathcal{G}_{E_T})=\bigcup_{n\geq0}\mathrm{Sym}(n+1).$$
\end{example}

We finish this section by considering briefly the case when the graph is acyclic. We have seen in Section \ref{sec:Leavitt} that it happens precisely when the Leavitt path algebra is finite dimensional and in our next result we state the topological full version of that fact: the graph is acyclic precisely when the associated topological full group is finite. This is well known to experts, but we have not been able to find it explicitly in the literature.

\begin{prop}\label{prop:acyclicfinite}
  Let $E$ be a finite graph and let $\Lambda_w$ denote the set of paths ending at the vertex $w$. Then, the following are equivalent:
  \begin{enumerate}
      \item[i)] $E$ is acyclic,
      \item[ii)] $\mathbf{F}(\mathcal{G}_{E})\cong\prod\limits_{w\in\operatorname{Sink}(E)}\operatorname{Sym}\left(\Lambda_w\right)$,
      \item[iii)] $\TFG$ is finite.
  \end{enumerate}
 where $\mathbf{F}(\mathcal{G}_{E})$ denotes the topological full group of the graph groupoid $\mathcal{G}_E$ and $\operatorname{Sym}(\Lambda_w)$ is the symmetric group on the paths ending at the vertex $w$.
\end{prop}
\begin{proof}
We first show that i) implies ii) (this is the topological full group version of \cite[Corollary 4.2.13]{AAS}). As $E$ is finite and acyclic it  has only finitely many paths, so its boundary set $\partial E=E^{\infty}\cup\left\{\mu\in E^{*}\mid r(\mu)\in E^{0}_{\text{sing}}\right\}$ is $$
\partial E=\bigsqcup\limits_{w\in\operatorname{Sink}(E)}\Lambda_w.
$$
and thus $\partial E$ is a finite set. This also implies that $\partial E$ is the maximal forest-basis in the sense that every forest-basis can be expanded to $\partial E$ and every element in the group $\TFG$ is a permutation of $\partial E$. The permutation has one orbit for each sink, and $|\Lambda_w|$ is the size of the orbit so we deduce that $\TFG$ can be decomposed as in ii).

As the groups in ii) are finite, to finish the proof we only have to show that iii) implies i). We argue by contradiction and assume that $E$ contains a cycle $C$ corresponding to a word $\alpha=\alpha_0\alpha_1\ldots\alpha_s\in E^f$.  We distinguish two cases, according to whether the cycle has an exit or not.

Assume first that $C$ has an exit $\beta$, i.e., there is $\beta\in E$ with $v=s(\beta)\in C$, $u=r(\beta)\not\in C$ and we may assume $v=s(\alpha_0)$.  There is a family of forest-bases  $B_n$, $n\geq 0$, so that for each $n$, $B_n$ is the smallest forest-basis containing $\alpha^n$. To construct the family, observe that $B_0=E^0$ is just the trivial forest-basis and to form $B_1$, we expand $u\in B_0$, then expand successively  $\alpha_0$, $\alpha_0\alpha_1$ and so on until we get a forest-basis containing $\alpha$, this basis is $B_1$ and we may repeat. In this process, when we pass from $B_n$ to $B_{n+1}$, we are adding one more path ending in $\beta$ and therefore having range $u$. This means that in each $B_n$ we have at least $n+1$ paths with range $u$. And therefore we have a copy of the symmetric group $\mathrm{Sym}(n+1)$ inside $\TFG$ consisting on elements of the form $(B_n,B_n,\varphi)$ for $\varphi$ a permutation between those paths. This contradicts iii).

  Now, we are left with the case when the cycle $C$ does not have an exit. Let $v$ be any vertex in $C$.
 Again, let $B_0=E^0$ be the trivial forest-basis and $B_1$ the smallest forest-basis containing $\alpha$. As $C$ has no exit, we see that $|B_0|=|B_1|$ and we have an element $g\in\TFG$ with forest-bases diagram $(B_1,B_0,\varphi)$ where $\varphi$ is the identity in $E^0\setminus\{v\}$ and sends $\alpha\mapsto v$. We claim that $g$ has infinite order, to see note that for $z=l(\alpha)$
$$
U=Z(\alpha,v)\sqcup \left(\mathcal{G}_E^{(0)}\backslash Z(v,v)\right)=\{(\alpha x,z,x)\mid x\in\partial E,s(x)=v\}\sqcup \left(\mathcal{G}_E^{(0)}\backslash Z(v,v)\right)
$$
and for any $k\in\Z$,
$$
U^k=\{(\alpha x,kz,x)\mid x\in\partial E,s(x)=v\}\sqcup \left(\mathcal{G}_E^{(0)}\backslash Z(v,v)\right)
$$
and the claim follows. Again, this contradicts iii).\end{proof}


\begin{remark} Note that it is crucial in this result that we are working with Nekrashevich's version of topological full groups.  If the graph does not have condition (L), i.e., when the groupoid is not effective (see Remark \ref{rem:L}), Matui's version could also be finite even if the graph is not acyclic. 
\end{remark}

\subsection{Embedding the graph topological full group in the Leavitt path algebra}

Next, we want to embed the topological full group $\TFG$ into the Leavitt path algebra $L_K(E)$. The fact that there exists such an embedding in the case of the  Higman-Thompson groups $V_{d,n}$  was first observed by Pardo \cite{Pardo}, see also \cite{DicksConch}.  For the general case, it is convenient to perform the embedding using Steinberg algebras that we define next.

Let $\mathcal G$ be a Hausdorff ample groupoid and $K$ a field. The \emph{Steinberg algebra} $A_K(\mathcal G)$ is the $K$-vector space of locally constant, compactly supported functions $f:\mathcal G\to K$, with pointwise addition and scalar multiplication, and convolution product
$$
(f*g)(\gamma)=\sum_{\substack{(\eta,\zeta)\in\mathcal G^{(2)}\\
\eta\zeta=\gamma}}
f(\eta)g(\zeta).
$$
Equivalently,
$$
A_K(\mathcal G)=\operatorname{span}_K\{\mathbf{1}_B:B\subseteq\mathcal G\text{ is a compact open bisection}\}
$$
where $\mathbf{1}_B$ is the indicator function taking value $1$ for elements in $B$, 0 otherwise. Steinberg algebras are sometimes considered generalized Leavitt algebras due to the fact that in the case of a graph groupoid $\mathcal G_E$ the Steinberg algebra $A_K(\mathcal G_E)$ and the Leavitt path algebra $L_K(E)$ are isomorphic (see \cite[Example 3.2]{OS}). This isomorphism can be stated explicitly as 
$$
\begin{aligned}
\Phi:L_K(E)&\longrightarrow A_K(\mathcal G_E),\\
\mu\nu^* &\longmapsto \mathbf1_{Z(\mu,\nu)}\qquad\text{whenever }r(\mu)=r(\nu).
\end{aligned}
$$
On the other hand, there is a natural representation of topological full groups of nice enough groupoids into the Steinberg algebra, sending a bisection $U$ to the function $\mathbf{1}_U$. 
 In the case we are interested in, i.e., groupoids $\mathcal{G}_E$ where the graph $E$ is finite we may realize this representation via the Leavitt path algebra as follows. Given $U\in\TFG$ and a forest-bases representative $(B_1,B_2,\varphi)$, put $B_1=\{\mu_1,\ldots,\mu_k\}$, $B_2=\{\nu_1,\ldots,\nu_k\}$ and $\varphi(\mu_i)=\nu_i$. Then set
$$\iota(U)=\sum_{i=1}^k\mu_i\nu_i^*\in L_K(E).$$
Obviously, the composition $\Phi\circ\iota$ is the natural representation of $\TFG$ in the Steinberg algebra and this implies
 

\begin{theorem}\label{thm:embedding} The map $\iota:\TFG\to L_K(E)$ is a well defined group monomorphism.
\end{theorem}

As observed above, there is an action of the topological full group $\TFG$ on $\partial E$ given for $U\in\TFG$ with $s,r:U\to\partial E$  by $\pi_U=s\circ r^{-1}$.
The element $\iota(U)$ defined above  allows one to make this action explicit because for any $x\in\partial E$,
$$\pi_U(x)=(\iota(U))x$$
where we operate in the word in the right hand side according to the rules of the Leavitt path algebra $L_K(E)$. We illustrate this in our next Example.

\begin{example}\label{ex:x0-action}
Let $E$ be the $2$-rose, with loops $y_1$ and $y_2$. Thus $\partial E=\{y_1,y_2\}^{\mathbb N}$ is the Cantor set. Consider the full bisection
$$
U=Z(y_1y_1,y_1)\sqcup Z(y_1y_2,y_2y_1)\sqcup Z(y_2,y_2y_2)
$$
that corresponds to the element $x_0$ in Figure (\ref{fig:generators-F}) below. We have
$$
\iota(U)=y_1y_1y_1^*+y_1y_2(y_2y_1)^*+y_2(y_2y_2)^*=y_1y_1y_1^*+y_1y_2y_1^*y_2^*+y_2y_2^*y_2^*.
$$
For every $z\in\partial E$, the induced homeomorphism is
$$
\pi_U(z)=(y_1y_1y_1^*+y_1y_2y_1^*y_2^*+y_2y_2^*y_2^*)z
$$
where depending on whether $z$ starts by $y_1$, $y_2y_1$ or $y_2y_2$ all the terms except of the first, the second, or the third, vanish.
\end{example}



\subsection{Filtration of the topological full group}

The embedding $\iota$ of Theorem \ref{thm:embedding}  induces the following filtration of the topological full group
$$\TFG_n=\{U\in\TFG\mid\iota(U)\in (L_K(E))_n\}.$$
This means that if a bisection $U\in\TFG$ admits some forest-bases expression $(B_1,B_2,\varphi)$ with $B_1=\{\mu_1,\ldots,\mu_k\}$, $B_2=\{\nu_1,\ldots,\nu_k\}$ and $\varphi(\mu_i)=\nu_i$ such that
$$\max_{i=1}^k(l(\mu_i)+l(\nu_i))\leq n,$$
then $U$ lies in $(\TFG)_n$. A priori, it could happen that this inequality fails but it is possible to modify $\iota(U)$ using Cuntz-Krieger relations so that we get an element in  $(\TFG)_n$. We show first that this can not be the case.


\begin{lemma}\label{lemma:SteinbergLeavittFilt}
Let $U\in\TFG$ and $n\geq0$. Then $U\in\TFG_n$ if and only if $U$ has a forest-bases representative $(B_1,B_2,\varphi)$ satisfying
$$
l(\mu)+l(\varphi(\mu))\leq n\qquad\text{for every }\mu\in B_1.
$$

\end{lemma}
\begin{proof} The if part is obvious, so we only have to show that this condition is necessary. We assume that  $\iota(U)$ can be transformed using the CK-relations into an expression lying in $L_K(E)_n$. As $\iota(U)$ has integer coefficients and the CK-relations have also integer coefficients, we see that  we can rewrite $\iota(U)$ as a finite linear combination with integer coefficients of monomials of the form $\alpha_j\beta_j^*$ such that $r(\alpha)=r(\beta)$ and $l(\alpha)+l(\beta)\leq n$,
\begin{equation}\label{eq:iotaU}
\iota(U)=\sum_{j=1}^{r}a_j\beta_j\alpha_j^*.
\end{equation}
We are going to modify this expression by transforming these monomials into others which are as long as possible. Consider one of them, say $\beta\alpha^*$. If $r(\beta)$ is a sink or if $l(\beta)+l(\alpha)\geq n-1$, we do nothing. In other case, we use the Cuntz-Krieger relation
\begin{equation}\label{eq:suma1}
\beta\alpha^*=\sum_{\substack{e\in E^1\\s(e)=r(\beta)}}(\beta e)(\alpha e)^*.
\end{equation}
We substitute this in (\ref{eq:iotaU}) and repeat the process until there is no monomial left for which this applies. Grouping together possible repeated monomials, this means that up to these modifications and possibly some  change of coefficients, we get an expression as (\ref{eq:iotaU}) such that all the monomials are pairwise distinct and for all of them we have that either $r(\beta_j)$ is a sink, or the total length $l(\beta_j)+l(\alpha_j)$ is at least $n-1$ (and in any case, it is at most $n$).

We claim that the bisections $Z(\beta_j,\alpha_j)$ for  $j=1,\cdots,r$ are pairwise disjoint. Indeed, in other case, one would be inside the other, meaning that the defining paths of the bigger one would be prefixes of the defining paths of the smaller one. But as our bisections are pairwise distinct, the condition on the lengths that we are assuming make this impossible so the claim follows.

Now, apply the map $\Phi:L_{K}(E)\to A_{K}(G_E)$. We get  
\begin{equation}\label{eq:suma4}
\mathbf1_U=\sum_{j=1}^{r}a_j\mathbf1_{Z(\beta_j,\alpha_j)}.
\end{equation}
But since $Z(\beta_j,\alpha_j)$ are pairwise disjoint and $\mathbf1_{U}, \mathbf1_{Z(\mu_i,\nu_i)}$ only take values $0$, $1$, we deduce that  the coefficients $a_{j}$ can only be 0 or 1, so we may assume removing some zeros if necessary that they are all equal to one. So we get
$$
U=\sqcup_{j=1}^{r}Z(\beta_j,\alpha_j)\,\,\text{ with }\,\,l(\beta_j)+l(\alpha_j)\leq n
$$
and since $U$ is a bisection, we deduce that $B_1=\{\beta_1,\ldots,\beta_r\}$ and $B_2=\{\alpha_1,\ldots,\alpha_r\}$ are forest-bases and we have the desired forest-bases representative.
\end{proof}

In our next result, we describe a general upper bound for $|\TFG_n|$.  This bound can be improved in many particular cases, we illustrate this with the case of the Toeplitz algebra and in Section \ref{sec:preThompson} we do the same for Thompson's group $V$.

\begin{theorem}\label{theorem:TFGbound} Let $E$ be a finite graph and denote by $f(n)$ the number of paths in $E^f$ of length at most $n$. Then
for every $n\ge0$,
$$
|\TFG_n|\leq\sum_{k=0}^{f(n)}k!\binom{f(n)}{k}^{\!2}\leq f(n)!\,2^{2f(n)}.
$$
\end{theorem}
\begin{proof}

By Lemma~\ref{lemma:SteinbergLeavittFilt}, we have that every $U\in(\TFG)_n$ has a forest-bases representative $(B_1,B_2,\varphi)$ such that for any $\alpha\in B_1$, $l(\alpha)+l(\varphi(\alpha))\leq n$ so in particular, $l(\alpha),l(\varphi(\alpha))\leq n.$ For each possible size of $B_1,B_2$, say $k$, there are $\binom{f(n)}{k}^2$ possible ways to choose two subsets of the set of words of length at most $n$ and there are $k!$ possible bijections between them. Of course, not all these choices yield elements in the topological full group, and even if they do, some could be reducible, but all these choices do give the upper bound $|(\TFG)_n|\leq\sum_{k=0}^{f(n)} k!\binom {f(n)}k^2$, Finally, observe that
$$
\sum_{k=0}^{f(n)} k!\binom {f(n)}k^2\leq f(n)!\sum_{k=0}^{f(n)}\binom {f(n)}k^2=f(n)!\binom{2f(n)}{f(n)}\leq f(n)!\,2^{2f(n)}.
$$
\end{proof}

Observe that, with the notation of the proof of Theorem \ref{thm:growth} the number $f(n)$ is
$f(n)=\sum_{\alpha=0}^{n}\sum_{j=1}^{|E^0|}\tau_j(\alpha)$. In fact, it follows from the proof of Theorem \ref{thm:growth} that $f(n)$ grows either exponentially or polynomially if the spectral radius $\rho(A_E)$ of the adjacency matrix is 1. Equivalently, $f(n)$ is the dimension of the $n$-th level of the standard filtration of the path algebra on $E$ so by \cite{entropy1}, $$\limsup_{n\to\infty}\frac{\log(f(n)-f(n-1))}{n}=\log(\rho(A_E)).$$
This means that when $f(n)$ is strictly exponential, the bound of Theorem \ref{theorem:TFGbound} grows much faster than in the case when $f(n)$ is polynomial. We illustrate this with two examples.


\begin{example}
Let $E=R_2$ be the 2-rose. As we have seen before, the associated topological full group is Thompson's group $V$. Let $(V)_n=V\cap(\L_2)_n$. Since $R_2$ has $2^j$ paths of length $j$, one has
$$
f(n)=\sum_{j=0}^n 2^j=2^{n+1}-1.
$$
and by the previous Theorem 
$$
|(V)_n|\leq\sum_{k=0}^{f(n)}k!\binom{f(n)}{k}^{\!2}\leq f(n)!\,2^{2f(n)}=(2^{n+1}-1)!\,2^{\,2(2^{n+1}-1)}.
$$
Now, consider the graph $E_T$ of Example~\ref{ex:Toeplitzgrowth} and let $A_{E_T}$ be its adjacency matrix. Then, $f(n)=2(n+1)$ for all $n\geq 0$ since $\sum_{i,j}(A^n_E)_{i,j}=2$ for all $n\geq 0$ because $A_{E_T}$ is idempotent. Applying Theorem~\ref{theorem:TFGbound}
$$
|\mathbf{F}(\mathcal{G}_{E_T})_n|\leq (2n+2)!2^{4(n+1)},\,\,\text{ for all }\,\, n\geq 0.
$$
\end{example}

In the next section we will consider again the group $V$ and prove that the bound above can be slightly improved.

\section{The filtration induced on Thompson's group \texorpdfstring{$V$}{V}.}
\label{sec:preThompson}

The simplest possible graph where the standard filtration of the Leavitt path algebra grows exponentially is the 2-rose $R_2$. In that case, as we have observed in Example \ref{ex:V} the associated topological full group $\mathbf{F}(\mathcal{G}_{R_2})$ is Thompson's group $V$. In this section we will look more closely at the filtration $(V)_n$ induced by the Leavitt path algebra embedding described in Section \ref{sec:TFG}.

\subsection{Tree pair diagrams and Thompson's groups \texorpdfstring{$V$}{V} and \texorpdfstring{$F$}{F}.}

We define the group $V$ as the topological full group of the 2-rose (see Example \ref{ex:Vgrowth}). As we have seen in Section \ref{sec:TFG}, elements in the topological full group of a finite graph are determined by a bijection between forest-bases. In the case of the 2-rose, forest-bases can be identified with finite rooted binary trees so it will be useful to fix some terminology. From now on, tree means  rooted binary tree. In such a tree,  nodes can have valence 3, 2 or 1. The only node of valence 2 is  the \emph{root}, nodes of valence 3 are called \emph{internal} and nodes of valence 1 are called \emph{leaves}. An internal node together with its $2$ outgoing edges and children is called a \emph{caret} and a caret whose children are
leaves is called an \emph{exposed caret}.

Any element of $V$ can be described as a triple $(T_1,T_2,\varphi)$ where $T_1$ and $T_2$ are trees and $\varphi$ is a bijection between their leaves. This representation is sometimes called \emph{tree pair diagram}. We will say that a tree pair diagram is \emph{reduced} if the corresponding forest-bases expression is, equivalently, if there is no other tree pair diagram representing the same element and having a smaller number of leaves.
 For example, the following reduced tree pair diagrams define elements $x_0,x_1$ in $V$ that will be important below.

\begin{figure}[h]
	\centering

	\begin{tikzpicture}[
			baseline=(current bounding box.center),
			every node/.style={font=\small},
			level distance=8mm,
			sibling distance=7mm
		]
		\node at (0,0.42) {$x_0$};

		\begin{scope}[xshift=-1.9cm,yshift=-0.1cm]
			\coordinate
			child {
					child {node {1}}
					child {node {2}}
				}
			child {node {3}};
		\end{scope}

		\draw[->,>=Stealth] (-0.65,-0.75) -- (0.65,-0.75);

		\begin{scope}[xshift=1.9cm,yshift=-0.1cm]
			\coordinate
			child {node {1}}
			child {
					child {node {2}}
					child {node {3}}
				};
		\end{scope}
	\end{tikzpicture}
	\qquad
	\begin{tikzpicture}[
			baseline=(current bounding box.center),
			every node/.style={font=\small},
			level distance=6mm,
			sibling distance=4mm
		]
		\node at (0,0.42) {$x_1$};

		\begin{scope}[xshift=-1.9cm,yshift=-0.1cm]
			\coordinate
			child {node {1}}
			child {
					child {
							child {node {2}}
							child {node {3}}
						}
					child {node {4}}
				};
		\end{scope}

		\draw[->,>=Stealth] (-0.65,-0.75) -- (0.65,-0.75);

		\begin{scope}[xshift=1.9cm,yshift=-0.1cm]
			\coordinate
			child {node {1}}
			child {
					child {node {2}}
					child {
							child {node {3}}
							child {node {4}}
						}
				};
		\end{scope}
	\end{tikzpicture}

	\caption{Tree pairs for $x_0$ and $x_1$.}
	\label{fig:generators-F}
\end{figure}
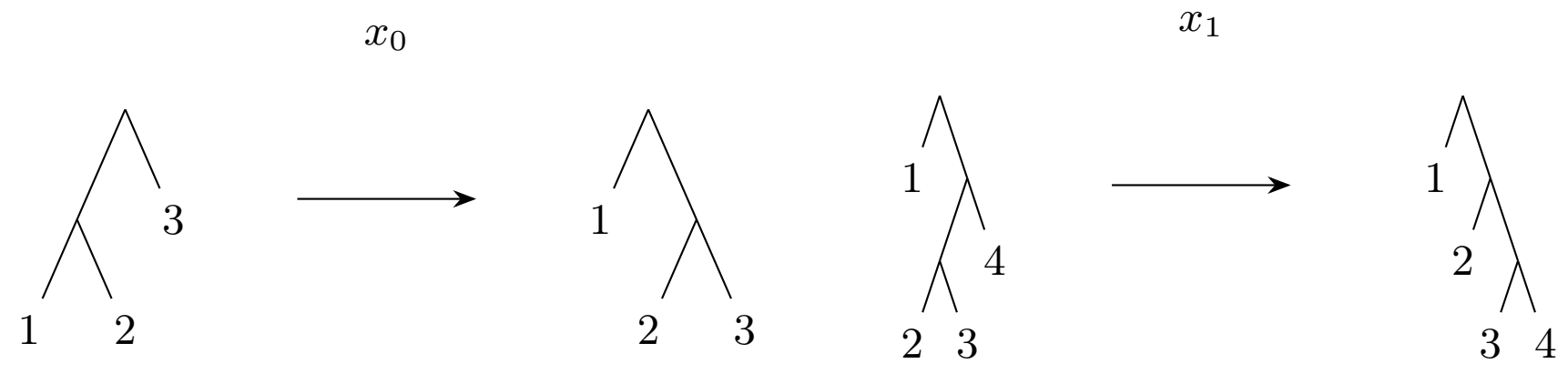

Observe that we may identify the set $R_2^f$ of path of finite length in the 2-rose  with the set of all leaves in trees, in what follows we will use this identification. In particular, if $a$ is a leaf in a tree, we will write $l(a)$ to denote its length as a path in $R_2$, equivalently, as a word in $y_1,y_2$. Therefore, the length $l(a)$ coincides with the distance between $a$ and the root in the usual tree metric.

The boundary $\partial R_2$ is the Cantor set $\mathfrak{C}$ so we get an action of the group $V$ on $\mathfrak{C}$.  In fact, the group $V$ is often defined in the literature using this action. 

The usual word length for elements in $V$ or $F$ is closely related to the number of carets in a reduced tree pair representation, see for example \cite{BurilloFelipe}. In our case,  we need to consider lengths of pairs of leaves. More precisely, 
using tree pair diagrams we get the following version of Lemma \ref{lemma:SteinbergLeavittFilt}

\begin{lemma}\label{lem:criV} Let $g\in V$ be an element with reduced tree pair diagram $(T_1,T_2,\varphi)$. Then $g\in(V)_n$ if and only if
$$\max_{a\text{ leaf of }T_1} \bigl(l(a)+l(\varphi(a))\bigr)\leq n. $$
\end{lemma}

The tree pair diagram representation can be useful to define the following two important objects: the subgroup $F$ of $V$ and the monoid $F^+$.

\begin{figure}[h]
\centering\begin{tikzpicture}[level 1/.style={sibling distance=7mm, level distance=7mm},
				level 2/.style={sibling distance=7mm, level distance=7mm},
				level 3/.style={sibling distance=7mm, level distance=7mm}
			]

			\coordinate
			child[draw=black] {node{}}
			child[draw=black] {
					child[draw=black] {node{}}
					child[draw=black] {
							child[draw=black] {node{}}
							child[draw=black] {node{}}
						}
				};
		\end{tikzpicture}
        \caption{A right-most tree.}
	\label{fig:rmtree}
\end{figure}
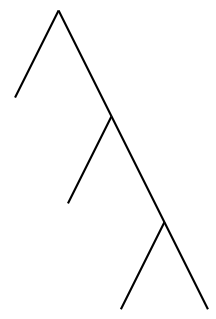

\begin{definition} Consider the lexicographic order induced by $y_1<y_2$ in the set of leaves of trees. The group $F$ consists of all the elements of $V$ for which the bijection in some (and therefore in any) tree pair diagram representation preserves the order, in that case we will denote the diagram just as $(T_1,T_2)$. The monoid $F^+$ consists of all the elements of $F$ that admit a tree pair diagram $(T_1,R)$ where $R$ is a \emph{right-most tree} (see Figure \ref{fig:rmtree}).
\end{definition}

For example, the elements $x_0$ and $x_1$ in \ref{fig:generators-F} lie in $F^+$ (so also in $F$).
In the literature, the group $F$ is often defined as a group of piecewise linear homeomorphisms of the unit interval on itself, and the study of the dynamics of this action has been very fruitful. In this paper we are more interested in combinatorial properties so the tree diagram representation of its elements will be more useful. We will also need the following well known fact.

\begin{theorem}\cite{CannonFloydParry1996} The elements $x_0$ and $x_1$ in Figure \ref{fig:generators-F} generate the group $F$. For $i\geq1$, set $x_i=x_0^{1-i}x_1x_0^{i-1}$. Then $F^+$ admits the following monoid presentation
$$F^{+}=\left\langle x_i\ (i\geq0)\ \middle|\ x_jx_i=x_ix_{j+1}\quad(i<j)\right\rangle_{\mathrm{mon}}.$$
\end{theorem}

In the next section we will also consider the following monoid 

\begin{definition} Let $M$ be the submonoid of $F^+$ generated by the elements $x_0$ and $x_1$.
\end{definition}

\subsection{Filtrations of \texorpdfstring{$F^+$}\ and \texorpdfstring{$M$}.}\label{subsec:monidfiltration}
In most cases, it is extremely difficult to determine exactly the number of elements in each term of the filtration induced by the standard filtration of the Leavitt path algebra in the associated topological full group. In this subsection, we use the fact that the tree pair diagram expression for elements of $F^+$ depends only on a single tree to show that things are  easier for $(F^+)_n$ and $(M)_n$. More precisely, we get a recursive formula for $|(F^+)_n|$ and a closed formula for  $|(M)_n|$.

We will use the following function. Let $T$ be a tree with $m$ leaves labeled  $a_1,\ldots,a_{m}$ in lexicographic order and set $\widehat{W}(T):=0$ if $T$ is the trivial tree (with just one node) and
\begin{equation}\label{eq:peso}
	\widehat{W}(T):=\max_{i=1,\dots,m} \left\{l(a_{i})+i\right\}
\end{equation}
otherwise.  Let $R_m$ be the right-most tree with $m$ leaves and $g\in F^{+}$ the element with  tree pair diagram $(T,R_m)$. This diagram is reduced if and only if the right-most caret of $T$ is not exposed, in that case we say that the tree $T$ is \emph{reduced} and this means that $l(a_m)<l(a_{m-1})$. 
Then, the lengths of the leaves of $R_m$ are $1,2,\ldots,m-1,m-1$ so the lengths of the monomials of $\iota(g)$ as in Lemma \ref{lem:criV} are
$$l(a_1)+1,l(a_2)+2\ldots,l(a_{m-1})+m-1,l(a_{m})+m-1.$$
But as we have $l(a_m)<l(a_{m-1})$, we deduce $l(a_m)+m\leq l(a_{m-1})+m-1$ thus
$$\widehat{W}(T)=\max\left\{l(a_1)+1,l(a_2)+2\ldots,l(a_{m-1})+m-1,l(a_{m})+m-1\right\}$$
so by Lemma \ref{lem:criV}, $g$ lies in $(F^+)_n$ if and only if $\widehat{W}(T)\leq n$.

Now, set
$$
	D_{m}^{n}:=\left|\left\{T\mid T\text{ tree with $m$ leaves, } \widehat{W}(T)\leq n\right\}\right|,
    S_{m}^{n}:=\left|\left\{T\mid T\text{ reduced tree with $m$ leaves, } \widehat{W}(T)\leq n\right\}\right|,
$$
in particular $D_1^n=S_1^n=1$, $S_2^n=0$ for $n\geq 1$ and $D_m^1=D_m^2=S^1_m=S^2_m=0$ for $m\geq 2$. 
\begin{prop}\label{prop:recurrencia}
    For $n\geq 3$, $|(F^+)_n|=\sum_{m=1}^\infty S_{m}^{n}$ and
    $$D^n_m=\sum_{i=1}^{m-1}D_i^{n-1}D_{m-i}^{n-1-i}\text{ for }m\geq 2,$$
    $$S^n_m=\sum_{i=1}^{m-1}D_i^{n-1}S_{m-i}^{n-1-i}\text{ for }m\geq 3.$$
\end{prop}
\begin{proof}  The first equation follows from  the previous considerations and the definition of $S_m^n$. 
Observe that in fact the sum is finite, because $S_{m}^{n}=0$ if $m>>n$.
For the other two equations, we can argue using a similar approach as one does for the Catalan numbers (see Subsection \ref{sec:double}). Specifically, to form all the non trivial trees $T$ with $m$ leaves and $\widehat{W}(T)\leq n$, start with a single caret and attach subtrees $T_1$ and $T_2$ to the left and right hand leaves respectively so that $T_1$ has $i$ leaves and $\widehat{W}(T_1)\leq n-1$, and $T_2$ has $m-i$ leaves and $\widehat{W}(T_2)\leq n-1-i$. This yields the first equation. For the second, proceed in the same way but take into account that for $T$ to be reduced, $T_2$ has to be reduced too.
\end{proof}

These formulas are not very informative so now we will consider the finitely generated submonoid $M$ instead. We begin with an easy observation.

\begin{lemma}\label{lem:Fplus-expansions}
	Let $1\neq g\in F^+$ be an element with reduced tree pair $(T,R_m)$, where $R_m$ denotes the right-most tree with $m$ leaves and
	$m\geq3$. Let $\Psi_i(T)$ be the tree obtained by attaching one caret to the $i$-th leaf of $T$. Then
	\begin{equation}\label{eq:Fplus-expansions}
		g\circ x_0=(\Psi_1(T),R_{m+1}),\qquad
		g\circ x_1=(\Psi_2(T),R_{m+1}),
	\end{equation}
	and both pairs are reduced.
\end{lemma}
\begin{proof}
	Right composition with $x_{0}$ or $x_{1}$ attaches a caret to
	the first or second leaf of $T$, respectively, and replaces $R_m$ with
	$R_{m+1}$. Since $T$ has at least 3 leaves, the resulting trees are reduced.
\end{proof}

Let $T_0$ and $T_1$ denote the source trees in $x_0$ and $x_1$, respectively.

\begin{prop}\label{prop:Fplus-source-recursion}
	Let $\Lambda_m$ be the set of all reduced trees $T$ with $m$ leaves such that $(T,R_m)$ lies in
$M$. Then $\Lambda_3=\{T_0\}$, $\Lambda_4=\{T_1\}\sqcup \{\Psi_1(T_0)\}\sqcup \{\Psi_2(T_0)\}$ and for $m\geq 4$,
	\begin{equation}\label{eq:Fplus-source-recursion}
					\Lambda_{m+1} =\Psi_1(\Lambda_m)\sqcup \Psi_2(\Lambda_m).
	\end{equation}
	Thus $|\Lambda_3|=1$ and $|\Lambda_m|=3\cdot2^{m-4}$ for $m\geq4$.
\end{prop}
\begin{proof} The assertion on $\Lambda_3$ and $\Lambda_4$ is obvious. For the inductive step, let $g\in M$ with reduced pair $(T,R_{m+1})$. We can express $g$ as a word in the monoid generators $x_0$ and $x_1$ so we let $g=hx$ where $x\in\{x_0,x_1\}$. This means that if $h$ has reduced pair $(\hat{T},R_m)$, then Lemma~\ref{lem:Fplus-expansions} implies that $T=\Psi_1(\hat{T})$ if $x=x_0$ and
$T=\Psi_2(\hat{T})$ otherwise. Observe that trees in the image of $\Psi_1$ have the left-most caret exposed whereas trees in the image of $\Psi_2$ have a left-most single leaf which does not belong to an exposed caret. This implies that both images are disjoint sets and also that the maps $\Psi_1$ and $\Psi_2$ are injective, because there is always an exposed caret whose contraction yields the original tree. Therefore $|\Lambda_{m+1}|=|\Psi_1(\Lambda_m)|+|\Psi_2(\Lambda_m)|=2|\Lambda_{m}|$
\end{proof}

\begin{theorem}\label{thm:Fplus-filtered-count}
	We have $1=\left|(M)_0\right|=\left|(M)_1\right|=\left|(M)_2\right|=\left|(M)_3\right|$, $2=\left|(M)_4\right|$ and
	\begin{equation}\label{eq:Fplus-filtered-count}
		\left|(M)_n\right|=3\cdot2^{n-5},\text{ for } n\geq5.
	\end{equation}
\end{theorem}
\begin{proof} We denote by $L_m$ the left-most tree with $m$ leaves. We start proving the following claim.

\noindent{\bf Claim:} Let $m\geq 3$. Then $\widehat{W}(L_m)=m+1$ and for any tree $T$ in $\Lambda_m$ such that $T\neq L_m$, we have $\widehat{W}(T)=m+2$.

The fact that $\widehat{W}(L_m)=m+1$ is obvious because the lengths of the leaves of $L_{m}$ are $m-1,m-1,m-2,\ldots,1$ so 
\begin{equation}\label{eq:leftmost}\widehat{W}(L_{m})=\max\{m-1+1,m-1+2,m-2+3,\ldots,1+m\}=m+1.\end{equation}
For the second assertion, argue by induction on $m$. The base case $m=3$ is vacuous, so the first non trivial base case is $m=4$. For $m=4$ we have that $\widehat{W}(L_4)=5$ and $\widehat{W}(T_1)=\widehat{W}(\Psi_2(L_3))=6$. Let $m\geq 5$ and $T\in\Lambda_m\backslash\{L_m\}$. By Proposition \ref{prop:Fplus-source-recursion} we can write $T=\Psi_j(\tilde{T})$ where $j$ is 1 or 2 and $\tilde{T}\in\Lambda_{m-1}$. 

Assume first that $\tilde{T}= L_{m-1}$. Then we must have $T=\Psi_2(\tilde{T})$ because in other case $T$ would be $L_m$. Therefore the lengths of the leaves of $\tilde{T}$ ordered lexicographically are
$$\{m-2,m-2,m-3\ldots,1\}$$
so the lengths of the leaves of $T$ are
$$\{m-2,m-1,m-1,m-3,\ldots,1\}$$
and $\widehat{W}(T)=\max\{m-2+1,m-1+2,m-1+3,m-3+4,\ldots,1+m\}=m+2$.

So we assume now that $\tilde{T}\neq L_{m-1}$.
Let $a_1,\ldots,a_{m-1}$ be the leaves of $\tilde{T}$ ordered lexicographically. By induction hypothesis we have
$$\widehat{W}(\tilde{T})=\max\{l(a_1)+1,l(a_2)+2,\ldots,l(a_{m-1})+m-1\}=m-1+2=m+1.$$
Assume first that $j=1$, i.e., $T=\Psi_1(\tilde{T})$. Then, the lengths of the leaves of $T$ are $$\{l(a_1)+1,l(a_1)+1,l(a_2),\ldots,l(a_{m-1})\}$$ thus
$$\widehat{W}(T)=\max\{l(a_1)+2,l(a_1)+3,l(a_2)+3,\ldots,l(a_{m-1})+m\}=1+\widehat{W}(\tilde{T})$$
where the last equality follows from the fact that for any tree, $l(a_1)\leq l(a_2)$ so the maximum is never attained at the second leaf. As by induction we are assuming $\widehat{W}(\tilde{T})=m+1$, we are done.

We are left with the case when  $T=\Psi_2(\tilde{T})$. Then, 
the lengths of the leaves of $T$ are $$\{l(a_1),l(a_2)+1,l(a_2)+1,l(a_3),\ldots,l(a_{m-1})\}$$ and so
$$
\begin{aligned}
\widehat{W}(T)&=\max\{l(a_1)+1,l(a_2)+3,l(a_2)+4,l(a_3)+4\ldots,l(a_{m-1})+m\}\\
&=\max\{\widehat W(\widetilde T)+1,\ l(a_2)+4\}.
\end{aligned}
$$
Since $\tilde{T}$ has $m-1$ leaves, $l(a_2)\leq m-2$. Thus 
$$
 \widehat W(T)=\max\{m+2,\ l(a_2)+4\}=m+2. 
$$

At this point, observe that the claim above implies that if $T$ belongs to $\Lambda_m$ (for $m\geq 3$), then $g=(T,R_m)\in(M)_n$ if and only if either $m\leq n-2$ or $T=L_{n-1}$. And as a consequence,

$$|(M)_n|=2+\sum_{m=3}^{n-2}|\Lambda_m|, \text{ for } n\geq 4$$
where the term 2 comes from the trivial element (that corresponds to the trivial tree) and the tree $L_{n-1}$. Using 
Proposition \ref{prop:Fplus-source-recursion} we get for $n\geq 5$
$$|(M)_n|=2+|\Lambda_3|+\sum_{m=4}^{n-2}|\Lambda_m|=3+3\sum_{m=4}^{n-2}2^{m-4}=3(1+\sum_{m=0}^{n-6}2^m)=3\cdot 2^{n-5}.$$
\end{proof}

\subsection{Double-exponential upper and lower bounds for \texorpdfstring{$|(V)_n|$}.}\label{sec:double}

In this subsection, we show Theorem \ref{teoInt:doubleexp}, i.e., we show that the $n$-th level of the filtration in $V$ induced by the standard filtration of the Leavitt algebra $\mathbb L_2$ admits double exponential lower and upper bounds.

To do that, we will introduce certain combinatorial numbers related to the numbers $S_k^n$ and $D_k^n$ of Subsection \ref{subsec:monidfiltration}, but this time it will be more convenient to index in terms of carets instead of leaves. So we set
$$C_k^n=|\{T\mid\text{ tree with $k$ carets and  leaves of length  at most }n\}|.$$
Recall that the Catalan numbers are defined as
$$C_k=|\{T\mid\text{ tree with $k$ carets}\}|.$$
As the leaves of a tree with $k$ carets can have length at most $k$,  $C_k=C_k^n$ if $n\geq k$. For $k\geq 1$, the Catalan numbers satisfy  the well known recurrence
$$C_k=\sum_{i=0}^{k-1}C_iC_{k-1-i}$$
which follows by observing that any tree with $k$ carets can be obtained starting with a single caret and attaching trees $T_1$ and $T_2$ so that $T_1$ has $i$ carets and $T_2$ has $j$ where $i+j=k-1$. With a similar argument, we get

\begin{prop}\label{prop:Ckn}
	Let $n\geq0$. Then, the numbers $C_k^n$ satisfy
	$C_0^n=1,\,\, C_k^0=0\,\,(k\geq1),$ $C_k^n>0$ whenever $0\leq k\leq2^n-1$ and
	\begin{equation}\label{eq:Ckn-rec}
		C_k^{n+1}=\sum_{i=0}^{k-1}C_i^nC_{k-1-i}^n
		\qquad(k\geq1).
	\end{equation}
	\end{prop}
\begin{proof} The fact that $C_0^n=1,\,\, C_k^0=0\quad(k\geq1),$ is obvious. A binary tree  whose leaves have length at most $n$ can have at most $2^n-1$ carets. Every number of carets between $1$ and $2^n-1$ can occur so we get the second assertion. Finally, for the recurrence, observe that we can proceed as for the Catalan numbers: to form any tree with $k$ carets and length of the leaves bounded by $n+1$, start with a single caret and attach trees $T_1$ and $T_2$ with $i$ and $j$ carets respectively for $i+j=k-1$ whose leaves have length at most $n$. This proves~\eqref{eq:Ckn-rec}. \end{proof}

Fixing $n$, the generating polynomial of the series $C_k^n$ has the following property.

\begin{prop}\label{prop:Ckn-generating-polynomial}
	The polynomials $F_n(t):=\sum_{k=0}^\infty C_k^n t^k =\sum_{k=0}^{2^n-1}C_k^n t^k$ satisfy
	\begin{equation}\label{eq:Fn-rec}
		F_0(t)=1,\qquad F_{n+1}(t)=1+tF_n(t)^2.
	\end{equation}
	As a consequence, 
	\begin{equation}\label{eq:Fn-nested}
		F_n(t)=1+t\left(1+t\left(1+t\left(\cdots(1+t)^2\cdots\right)^2\right)^2\right)^2.
	\end{equation}
\end{prop}
\begin{proof} We have
	$$
		\begin{aligned}
			\left(F_{n}(t)\right)^{2} =\left(\sum_{k=0}^\infty C_{k}^{n}t^{k}\right)^{2} =\sum_{k=0}^\infty\left(\sum_{i=0}^{k}C_{i}^{n}C_{k-i}^{n}\right)t^{k}  =\sum_{k=0}^\infty C_{k+1}^{n+1}t^{k}                              =\frac{1}{t}\sum_{k=0}^\infty C_{k+1}^{n+1}t^{k+1}  &                    \\
			 =\frac{1}{t}\sum_{k=1}^\infty C_{k}^{n+1}t^{k}=\frac{F_{n+1}(t)-C_{0}^{n+1}}{t}                                =\frac{F_{n+1}(t)-1}{t}.&\\
		\end{aligned}
	$$
	Therefore $F_{n+1}(t)=1+t(F_n(t))^2$, and iterating this identity gives the nested
	expression~\eqref{eq:Fn-nested}.
\end{proof}

For instance,
$$
	F_3(t)=1+t\bigl(1+t(1+t)^2\bigr)^2 =1+t+2t^2+5t^3+6t^4+6t^5+4t^6+t^7.
$$

Set also
$$
	x_n:=\sum_{k=0}^\infty C_k^n,
$$
i.e., $x_n$ is the number of trees with leaves of length at most $n$. Obviously, $x_0=1$ and using again a similar argument as for the Catalan numbers and seeing a non trivial tree as the result of attaching two trees to a single caret yields
\begin{equation}\label{eq:xn-rec}
	x_{n+1}=1+x_n^2.
\end{equation}
This series has been considered in \cite[Example 2.2]{Aho}, where it is shown that there is a constant
$\delta=1.502837\cdots$ such that
\begin{equation}\label{eq:xn-explicit}
	x_n=\left\lfloor\delta^{2^n}\right\rfloor.
\end{equation}

We can now show the following technical result from which we will deduce Theorem \ref{teoInt:doubleexp}.

\begin{prop}\label{prop:simetrico}
	For every $n\geq0$,
	\begin{equation}\label{eq:simetrica}
		\sum_{k\geq 0}\left(C_k^{\lfloor n/2\rfloor}\right)^2k!
		\leq \left|(V)_n\right|
		\leq \sum_{k=0}^{\lfloor2^{n/2}\rfloor-1}\left(C_k^n\right)^2(k+1)!.
	\end{equation}
\end{prop}
\begin{proof}
Observe first that by Lemma \ref{lem:criV}, if $g\in V$ is an element with reduced tree pair diagram $(T_1,T_2,\varphi)$, then $g\in(V)_n$ if and only if
\begin{equation}\label{eq:crit}
\max_{a\text{ leaf of }T_1} \bigl(l(a)+l(\varphi(a))\bigr)\leq n. 
\end{equation}
	Now, for the lower bound, assume that $T_1$ and $T_2$ are trees with $k$ carets such that all the leaves have length at most $\lfloor n/2\rfloor$. Then (\ref{eq:crit}) holds for every bijection $\varphi$ between their leaves but the diagram could be reduced meaning that when we let $k$ vary we could have repeated elements.   
	Fix a pair $T_1,T_2$ as before, observe that both trees have $L=k+1$ leaves. Let $s,r$ be the number of exposed carets in $T_1$ and $T_2$ respectively. Then fixing an exposed caret in each of the trees $T_1$, $T_2$, we may form at most $(k-1)!$ non reduced tree pair diagrams for those trees. Hence, for $k\geq1$, the number of reducible bijections is at most $s\,r\,(k-1)!$. Since $s,r\leq k$, there are at least
	$$	(k+1)!-s\,r\,(k-1)! \geq (k+1)!-k^2(k-1)!=k! $$
possible ways to choose the bijection $\varphi$ so that the tree pair diagram is	reduced. These choices yield	distinct elements in $V$, so we get the left-hand side inequality in~\eqref{eq:simetrica}.

For the upper bound, let $g$ be an element of $(V)_n$ whose reduced tree pair diagram is $(T_1,T_2,\varphi)$ where $\{u_1,\cdots,u_m\}$ and $\{\varphi(u_1),\cdots,\varphi(u_m)\}$ are the leaves of $T_1$ and $T_2$ respectively. Hence, each of the trees in the diagram has $m$ leaves. If we set $a_i:=l(u_i)$, $b_i:=l(\varphi(u_i))$ for $i=1,\cdots, m$, by 
(\ref{eq:crit}) we get $a_{i}+b_i\leq n$ for $i=1,\cdots, m$.

\noindent{\bf Claim:} A tree of a tree pair diagram of an element of $(V)_n$ has at most $\lfloor2^{n/2}\rfloor$ leaves. To see it, note first that the length of the leaves of a binary tree $T$ satisfy
$$
\sum_{u \text{ leaf of } T} 2^{-l(u)}=1.
$$
This can be easily seen by assigning weight 1 to the root and splitting the weight of each vertex into its two children. Applying this fact to $T_1$ and $T_2$ one gets
$$
\sum_{i=1}^{m}2^{-a_i}=1=\sum_{i=1}^{m}2^{-b_i}.
$$
Hence, by the Cauchy-Schwarz inequality
$$
m\leq \sum_{i=1}^{m}2^{\left(\frac{n-(a_i+b_i)}{2}\right)}\leq 2^{n/2}\,\sum_{i=1}^{m}2^{-\left(\frac{a_i+b_i}{2}\right)}\leq 2^{n/2}\left(\sqrt{\sum_{i=1}^{m}2^{-a_i}}\right)\left(\sqrt{\sum_{i=1}^{m}2^{-b_i}}\right)=2^{n/2}.
$$
Since $m$ is an integer, we get that $m\leq\lfloor2^{n/2}\rfloor$, proving the claim.

Now, for each $k$ there are at most $(C_{k}^{n})^2$ choices for the pair of trees and at most $(k+1)!$ bijections between their leaves. Moreover the claim above implies that we only have to take into account tree pairs with at most $\lfloor2^{n/2}\rfloor$ leaves, i.e., at most $\lfloor2^{n/2}\rfloor-1$ carets.
Summing all the possibilities for $0\leq k\leq\lfloor2^{n/2}\rfloor-1$ we get the inequality in the right hand side.
\end{proof}

\bigskip

\noindent {\sl Proof of Theorem \ref{teoInt:doubleexp}.} For the lower bound, observe that (\ref{eq:xn-explicit}) implies
$$\left\lfloor\delta^{2^{\lfloor n/2\rfloor}}\right\rfloor=x_{\lfloor n/2\rfloor}=
\sum_{k\geq0}\left(C_k^{\lfloor n/2\rfloor}\right)\leq \sum_{k\geq0}\left(C_k^{\lfloor n/2\rfloor}\right)^2k!$$
so we only have to apply Proposition \ref{prop:simetrico}. For the upper bound, notice that by Proposition \ref{prop:simetrico} we get
	$$
\left|(V)_n\right|
		\leq\sum_{k=0}^{\lfloor 2^{n/2}\rfloor-1}\left(C_k^n\right)^2(k+1)!.
	$$
    Now, by the claim in the proof of Proposition \ref{prop:simetrico}, for the number $k$ of carets  we have $k\leq \lfloor 2^{n/2}\rfloor-1$ and thus
    $$
\sum_{k=0}^{\lfloor 2^{n/2}\rfloor-1}\left(C_k^n\right)^2(k+1)!\leq \left(\lfloor 2^{n/2}\rfloor\right)! \sum_{k=0}^{\lfloor 2^{n/2}\rfloor-1}\left(C_k^n\right)^2\leq \left(\lfloor 2^{n/2}\rfloor\right)! \,x_n^2. 
    $$
It only remains to prove that $\left(\lfloor 2^{n/2}\rfloor\right)!\leq x_n$.  Notice that since $\frac{n}{2}2^{n/2}\leq 2^{n-1}$, for $n\geq 4$ we have 
$$
\left(\lfloor 2^{n/2}\rfloor\right)!\leq \left(\lfloor 2^{n/2}\rfloor\right)^{\left(\lfloor 2^{n/2}\rfloor\right)}\leq 2^{\frac{n}{2}2^{n/2}}\leq 2^{2^{n-1}}\quad\text{for }\,n\geq 4.
$$
Finally, the recurrence $x_{n+1}=1+x_{n}^{2}$ with $x_1=2$ gives $x_{n}\geq 2^{2^{n-1}}$, yielding the upper bound for $n\geq 4$. For $n=0,1,2,3$ the result follows by direct computation:
$$
(x_0,x_1,x_2,x_3)=(1,2,5,26), \qquad \left(\left\lfloor 2^{0}\right\rfloor!,\left\lfloor 2^{1/2}\right\rfloor!,\left\lfloor 2^{2/2}\right\rfloor!,\left\lfloor 2^{3/2}\right\rfloor!\right)=(1,1,2,2).
$$

\begin{remark}\label{remark:F-double-exponential}
 The lower bound in Theorem  \ref{teoInt:doubleexp} also holds for $n\geq 6$ when $V$ is replaced by $F$:
 $$\left\lfloor\delta^{2^{\lfloor n/2\rfloor}}\right\rfloor
		\leq \left|(F)_n\right|
		\qquad(n\geq6),	$$
with the same constant $\delta$. To see it, set
$$y_n=\left|\{g=(T_1,T_2)\in F\mid\text{ reduced with lengths of the leaves in  $T_1$ and $T_2$ at most }n\}\right|.$$
 One can check that $y_3=35>x_3=26$. Moreover, given $g=(T_1,T_2),h=(L_1,L_2)\in F$ reduced so that $T_1,T_2,L_1,L_2$ all have lengths of leaves bounded by $n$, we may form a new reduced tree pair diagram that yields an element in $F$ where the lengths of the leaves are at most $n+1$ just attaching $T_1$ and $L_1$ to a single caret, doing the same for $T_2$ and $L_2$ and extending the bijections that determine $g,h$. This implies $y_{n+1}\geq y_n^2$. Since $x_{n+1}=x_n^2+1$, induction gives
	$$y_n\geq x_n+1\qquad(n\geq3).$$
 
Applying this for $\lfloor n/2\rfloor$ yields
	$$
		x_{\lfloor n/2\rfloor}\leq y_{\lfloor n/2\rfloor}\leq|(F)_n|.
	$$    
For $n\geq 6$, this together with (\ref{eq:xn-explicit}) prove the lower bound. By contrast, Theorem \ref{thm:Fplus-filtered-count} gives
$$
|M_{n}|=3\cdot 2^{n-5}\quad(n\geq5).
$$
Both $F$ and $M$ have the same generating system $S=\{x_0,x_1\}$, but the fact that $S$ generates $F$ as a group and $M$ only as a monoid has dramatic effects on the size of the corresponding filtrations.

\end{remark}






\end{document}